\documentclass[11pt]{article}
\usepackage{amssymb,amsmath,amsthm,mathrsfs,txfonts,color}
\def\pd#1#2{\frac{\partial#1}{\partial#2}}
\newtheorem{theorem}{\indent Theorem}[section]
\newtheorem{lemma}{\indent Lemma}[section]

\newtheorem{definition}{\indent Definition}[section]
\newtheorem{remark}{\indent Remark}[section]

\allowdisplaybreaks

\begin{document}
\title{Propagating Profiles of a Chemotaxis Model with Degenerate Diffusion:
Initial Shrinking, Eventual Smoothness and Expanding}
\author{
Tianyuan Xu$^a$, Shanming Ji$^b$\thanks{Corresponding author, email:jism@scut.edu.cn},
Ming Mei$^{c,d}$, Jingxue Yin$^a$,
\\
\\
{ \small \it $^a$School of Mathematical Sciences, South China Normal University}
\\
{ \small \it Guangzhou, Guangdong, 510631, P.~R.~China}
\\
{ \small \it $^b$School of Mathematics, South China University of Technology}
\\
{ \small \it Guangzhou, Guangdong, 510641, P.~R.~China}
\\
{ \small \it $^c$Department of Mathematics, Champlain College Saint-Lambert}
\\
{ \small \it Quebec,  J4P 3P2, Canada, and}
\\
{ \small \it $^d$Department of Mathematics and Statistics, McGill University}
\\
{ \small \it Montreal, Quebec,   H3A 2K6, Canada}
}
\date{}

\maketitle

\begin{abstract}
We investigate the propagating profiles of a degenerate chemotaxis model describing the
bacteria chemotaxis and consumption of oxygen by aerobic bacteria, in particular,
 the effect of the initial attractant distribution on bacterial clustering.
We prove that the compact support of solutions may shrink if the signal concentration satisfies a special structure,
and  show the finite speed propagating property without assuming the special structure on attractant concentration,
and obtain an explicit formula of the population spreading speed in
terms of model parameters.
The presented results suggest that bacterial cluster formation can be affected
by chemotactic attractants and density-dependent dispersal.


\end{abstract}

\section{Introduction}
We consider the following chemotaxis model with chemotactic consumption and porous media
diffusion
\begin{equation}\label{eq-modelnonboundary}
\begin{cases}
u_t=\nabla(\phi(u)\nabla u)-\chi \nabla\cdot(u\nabla v),\\
v_t=\Delta v-\alpha uv,\quad &x\in \Omega,~t>0,
\end{cases}
\end{equation}
where $u$ represents the number per unit volume of aerobic bacteria cells,
$v$ denotes the oxygen concentration, $\chi$ is the
chemotactic coefficient,
$\alpha$ denotes the fractional rate of oxygen consumption
per unit concentration of bacteria cell.
The diffusion of species is considered to
be degenerate in the form of $\nabla(\phi(u)\nabla u)$ with $\phi(u)=D u^{m-1}$
and $m>1$, which is dependent of the population density due to the population pressure.
This model can also describe other
chemotaxis progress with nutrient consumption.

In many biological cases, the diffusion coefficient $\phi(u)$ is not constant, which can be regarded as a consequence of the interaction between cells \cite{Bian2013Dynamic,Yao2013Blow,Mansour2008Traveling,Sherratt1996Nonsharp}.
It is worthy of mentioning that the porous medium type diffusion can represent
``population pressure'' in cell invasion models  \cite{Sherratt}, which initially arises from the ecology literature \cite{Gurney1975The,Gurtin1977On,Murry,XuMBE}.
In fact,
experimental investigation has shown that the diffusion coefficient depends on the bacterial density \cite{Wakita1994Experimental}.
In the bacterial experiments done by Ohgiwari {\it et al.} \cite{Ohgiwari1992Morphological},
they recognized that
cells
located inside the bacterial colonies move actively,
but cells became sluggish at the outermost front with
apparently low cell density.
This phenomenon indicates that bacteria become active as the cell density $u$ increases.
Thus,  a natural choice of the bacterial diffusion coefficient is $\phi(u)=u^{m-1}(m>1)$, and this porous medium type bacterial diffusivity is based on the degenerate  diffusion model proposed by Kawasaki {\it et al.} \cite{Kawasaki1997Modeling}.
Recently, Leyva {\it et al.} \cite{Leyva2013The}
incorporates a chemotactic term into the original model by Kawasaki {\it et al.}, and explores the effects of chemotaxis on bacterial aggregation patterns.

Chemotaxis
is the biased migration in the direction of a chemical
stimulus concentration gradient \cite{11,Hillen2009A}.
Bacteria can sense a large range of chemical
signals, such as the concentrations of nutrients, toxins, oxygen, minerals, etc.
A mathematical model for the process of aerobic motile bacteria toward oxygen which they consume was first proposed in \cite{Rosen1978Steady}. When $m=1$ in \eqref{eq-modelnonboundary}, namely,
the diffusion of bacteria cell are assumed to be random,
Tao and Winkler \cite{Tao2012Eventual} proved that this model admits a global weak solution, and a more interesting fact is that, the weak solution will become smooth after some time.
Recently,
chemotaxis models featuring a
density-dependent diffusion term have drawn great attention from many authors \cite{Francesco2008Fully,Bellomo2016A,Sugiyama2007Time,Wang2012Global,
Winkler2018Global,Manjun2015Global,Wang2007Classical}.
For this system with the porous medium diffusion (i.e. $m>1$ in \eqref{eq-modelnonboundary}),
it was shown that the weak solution is globally solvable in two dimension for $m>1$ \cite{Tao2013Global}.
In three dimensional space, the authors made great efforts to prove the
global existence of weak solutions for this model for any $m>1$.
Winkler and Tao \cite{Tao2013Locally,Winkler2015Boundedness} proved this problem admits a global weak solution for the case $m\in (1,\frac{8}{7}]$.

The interaction between diffusion and chemotaxis contributes substantial influences on the behavior of solutions for chemotaxis model with degenerate diffusion.
In \cite{Burger2006The}, Burger, Di Francesco, and Dolak considered the Keller-Segel model of chemotaxis with volume filling effect, which is degenerate when bacteria densities approaching either $0$ or $1$,
and they investigated the qualitative behavior of solutions, such as finite speed of propagation and asymptotic behavior of solutions.
Kim and Yao \cite{Kim2012The} studied the qualitative properties of the
Patlak-Keller-Segel model with porous medium type diffusion term by using
maximum principle type arguments,
and they proved the finite propagation property of the compactly supported solutions generated by
this type of degeneracy of diffusivity.
In \cite{Fischer2013Advection},
Fischer
proved finite speed of support propagation for the parabolic-elliptic chemotaxis Keller-Segel system with porous medium
type diffusive term and gave sufficient
criteria for support shrinking, based on the integral estimates and the Stampacchia's lemma.

The main feature of our model \eqref{eq-modelnonboundary} lies in that the porous medium diffusive term and the chemotactic term are in competition.
The dispersal term induces forward motion, whereas the chemotactic attraction may account for cohesive swarm and induce backward motion of the invasion boundary \cite{Fischer2013Advection}.
We explore the effect of density dependent diffusion and chemotactic attraction,
which can account for cohesive,
finite swarms with realistic density profiles.

To understand how changes in the initial conditions of chemotaxis can so dramatically alter the
aggregation behavior of bacteria, we study the effects of attractant concentration on
bacteria distribution. We will give (see
Theorem \ref{th-shrinking}) a mathematical understanding of the collective behavior of bacteria chemotactic toward oxygen.
We find that under certain initial conditions,
the boundary of $\text{supp}\, u(\cdot, t)$ moves backward in response to the gradient of attraction at early stage.
This indicates that
the size of the swarm is defined by a balance of chemotactic attraction and cell dispersal: the greater the attraction the smaller its size for a given total number of organisms.
This is observed biologically: bacteria
 exhaust the local oxygen and then
react to the attractant gradient they have created,
producing a flux towards the region
with more oxygen. Early in bioconvection, this process generated accumulations of cells, resulting in smaller size of cell collective region.
This experiment was conducted on  {\it Bacillus subtilis}
\cite{Cisneros2008The}.

One of the intrinsic characteristics of porous medium diffusions is the population moves with a finite speed of propagation,
which seems
more reasonable than infinite speed in biological applications. To put it concisely,
for any non-zero initial data $u_0$, the solution of linear diffusion equation
$u(x,t) > 0$ for $t > 0$ and any $x\in\mathbb R^N$,
thus a linear diffusion process predicts an
infinite propagation \cite{PME}.
However, the spatial support of the solution to the degenerate diffusion equation
remains bounded for all time $t > 0$ \cite{Fischer2013Advection}.

Bacteria are known to exhibit very diverse morphological aggregation patterns depending on a variety of environmental
conditions \cite{Ben1994Generic,Ohgiwari1992Morphological,Matsuyama1993Fractal,Wakita1994Experimental}.  These experimental observations showed the
bacterial envelop front propagate outward gradually over time
and the velocity of front propagation is finite.
In order to explain these phenomena, a variety of mathematical models have been proposed \cite{Kawasaki1997Modeling,Sherratt2010On,Satnoianu2012Travelling,Leyva2013The}.
The density-dependent degenerate diffusion model may
capture more pattern features found experimentally and
provides a better match to experimental cell density
profiles.
The difference between these
 diffusion types is that the porous medium type diffusion  leads to distinct boundaries, and the population density decreases
to zero at a finite point in space, rather than tends to zero asymptotically.
It is therefore not surprising that the behavioral property of living organisms in these two models
is different.
The porous medium type models allow the cells aggregate rather than spread out.
The non-physical
diffusion is eliminated in this model.

Although the underlying dynamics of the chemotaxis model with degenerate mobility can be complicated,
explicit description of bacteria invasion process can be given.
The challenge in the mathematical analysis consists of the chemotactic term as well as
the degeneracy of the diffusion term which generates compactly supported solutions.
We prove several propagating properties of solutions,
including the initial shrinking, finite propagation property, eventual smoothness and
eventual expanding.
The spreading
speed is the rate at which the species with uniformly positive
initial distribution over a large interval and zero distribution outside an interval
expands its spatial range \cite{Lewis2012Spreading}.
Theorem \ref{th-exact} below provides an explicit formula for the spreading speed in terms of model parameters.
To the best of our knowledge, it is the first work that presents a precise description of the propagating speed
for this model.
These results provide important insight into the spatial patterns and
rates of invading bacteria species in space.

The outline of this paper is as follows.
In Section 2, we state our main results and some notations.
We leave the comparison principle of the
corresponding degenerate chemotaxis equation
and its H\"older continuity into Section 3 as preliminaries.
Section 4 is devoted to the study of  the propagating properties of bacteria cells
and the large time behavior of the weak solution.

\section{Main results and Notations}

We consider the following chemotaxis system \eqref{eq-model} with degenerate diffusion
\begin{equation}\label{eq-model}
\begin{cases}
u_t=\Delta u^m-\chi\nabla\cdot(u\nabla v),\quad &x\in \Omega,~t>0,\\
v_t=\Delta v-uv,\quad &x\in \Omega,~t>0,\\[1mm]
\displaystyle
\pd{u^m}{n}=\pd{v}{n}=0,
\quad &x\in \partial\Omega,~t>0,\\
u(x,0)=u_0(x),\quad v(x,0)=v_0(x),\quad &x\in\Omega,
\end{cases}
\end{equation}
where $m>1$, $\chi>0$, $\Omega\subset\mathbb{R}^N$ is a
bounded domain with smooth boundary
and the spacial dimension $N\in\{1,2,3\}$,
$u_0,v_0$ are nonnegative functions,
$n$ is the unit outer normal vector.

Since degenerate diffusion equations may not have
classical solutions in general, we need to formulate the
following definition of generalized solutions
for the initial boundary value problem \eqref{eq-model}.

\begin{definition}
Let $T\in(0,\infty)$.
A pair of $(u,v)$ is said to be a weak solution
to the problem \eqref{eq-model} in $Q_T=\Omega\times(0,T)$ if

{\rm(1)} $u\in L^\infty(Q_T)$, $\nabla u^m\in L^2((0,T);L^2(\Omega))$,
and $u^{m-1}u_t\in L^2((0,T);L^2(\Omega))$;

{\rm(2)} $v\in L^\infty(Q_T)\cap L^2((0,T);W^{2,2}(\Omega))%
\cap W^{1,2}((0,T);L^2(\Omega))$;

{\rm(3)} the identities
\begin{eqnarray*}
&&\int_0^T\int_\Omega u\psi_tdxdt
+\int_\Omega u_0(x)\psi(x,0) dx
\\
&&\ \ =\int_0^T\int_\Omega \nabla u^m\cdot\nabla\psi dxdt
-\int_0^T\int_\Omega \chi u\nabla v\cdot\nabla\psi dxdt, \\
&& \int_0^T\int_\Omega  v_t\varphi dxdt
+\int_0^T\int_\Omega \nabla v\cdot\nabla\varphi dxdt
=\int_0^T\int_\Omega wz\varphi dxdt,
\end{eqnarray*}
hold for all $\psi,\varphi\in L^2((0,T);W^{1,2}(\Omega))%
\cap W^{1,2}((0,T);L^2(\Omega))$ with $\psi(x,T)=0$ for $x\in\Omega$;

{\rm(4)} $v$ takes the value $v_0$ in the sense of trace at $t=0$.

If $(u,v)$ is a weak solution of \eqref{eq-model} in $Q_T$ for
any $T\in(0,\infty)$, then we call it a global weak solution.

A pair of $(u,v)$ is said to be a globally bounded weak solution
to the problem \eqref{eq-model} if there exists a positive constant $C$ such that
\begin{equation*}
\sup_{t\in\mathbb R^+}\left\{\|u\|_{L^\infty(\Omega)}+\|v\|_{W^{1,\infty}(\Omega)}\right\}\le C.
\end{equation*}
\end{definition}

Throughout this paper we assume that the initial data satisfies
\begin{equation} \label{eq-initial}
u_0\in C(\overline\Omega), ~\nabla u_0^m\in L^2(\Omega), ~
\pd{v_0}{n}=0 \text{~on~}\partial\Omega, ~
v_0\in C^{2,\alpha_0}(\overline\Omega) \text{~for~some~} \alpha_0\in(0,1).
\end{equation}

We are aiming at the propagating properties of the cell invasions.
Let us first focus on the waiting time and the initial shrinking of the compact support
caused by chemotaxis.
Our approach is based on the comparison principle
and the technique of self-similar weak lower and upper solutions with compact support.

\begin{theorem}[Initial shrinking caused by chemotaxis] \label{th-shrinking}
Let $(u,v)$ be a globally bounded weak solution of \eqref{eq-model}
with (i) $N=1$; or
(ii) $\sup_{t\in(0,\infty)}\|u(\cdot,t)\|_{C^{1/(2m)}(\overline\Omega)}\le C$ for some
constant $C>0$.
Further we assume that
\begin{align} \label{eq-shrinking-1}
&\text{supp}\,u_0\subset\overline B_{R_0}(x_0)\subset\Omega,
\quad u_0\le K_0(R_0^2-|x-x_0|^2)^{d_0}, ~ x\in B_{R_0}(x_0),
\\ \label{eq-shrinking-2}
&\nabla v_0\cdot (x-x_0) \le -\mu |x-x_0|^2, ~ x\in B_{R_0}(x_0),
\end{align}
for some $x_0\in\Omega$ and positive constants $d_0\ge 1/(m-1)$, and $R_0,K_0,\mu>0$.
such that $\chi\mu>\frac{4m}{m-1}K_0^{m-1}\max\{1,R_0^{2((m-1)d_0-1)}\}$.
Then there exist a family of shrinking open sets
$\{A(t)\}_{t\in(0,t_0)}$ with $t_0>0$ such that $A(0)=B_{R_0}(x_0)$ and
$$\text{supp}\,u(\cdot,t)\subset\overline A(t)\subset\Omega, \quad t\in(0,t_0),$$
and $\partial A(t)$ has a finite negative derivative with respect to $t$.
\end{theorem}

\begin{remark}
The existence of globally bounded weak solutions of \eqref{eq-model}
is proved in \cite{JinarXiv}.
We will prove in Lemma \ref{le-L4} that (i) implies (ii).
The finite propagating speed (i.e. the derivative of $\partial A(t)$ with respect to $t$)
is interpreted as in the sense of Theorem \ref{th-exact}.
\end{remark}

We show the finite speed propagating property without the
special structure \eqref{eq-shrinking-2} on signal concentration.

\begin{theorem}[Finite speed propagating] \label{th-finite}
Let the assumptions in Theorem \ref{th-shrinking} be valid
except for \eqref{eq-shrinking-2}.
Then there exist a family of open sets
$\{A(t)\}_{t\in(0,t_0)}$ with $t_0>0$ such that $A(0)=B_{R_0}(x_0)$ and
$$\text{supp}\,u(\cdot,t)\subset\overline A(t)\subset\Omega, \quad t\in(0,t_0),$$
and $\partial A(t)$ has a finite derivative with respect to $t$.
\end{theorem}

\begin{remark}
Without the structure \eqref{eq-shrinking-2} on signal concentration,
we do not know the shrinking or expanding of the cells.
However, Theorem \ref{th-finite} shows the propagating speed is finite.
\end{remark}

If the cell density and the signal concentration have special structure,
we will present the exact propagating speed as follows.

\begin{theorem}[Exact propagating speed] \label{th-exact}
Let $(u,v)$ be a globally bounded weak solution of \eqref{eq-model}
with (i) $N=1$; or
(ii) $\sup_{t\in(0,\infty)}\|u(\cdot,t)\|_{C^{1/(2m)}(\overline\Omega)}\le C$ for some
constant $C>0$.
Further we assume that the initial values satisfy
\begin{equation} \label{eq-speed-i}
\begin{cases}
u_0= K_0\big[(R_0^2-|x-x_0|^2)_+\big]^{d}, ~ x\in \Omega,
\\
\nabla v_0\cdot (x-x_0)= -\mu |x-x_0|^2, ~ x\in B_{R_0}^\delta(x_0),
\end{cases}
\end{equation}
for some $x_0\in\Omega$ and positive constants $d=1/(m-1)$, $R_0,K_0,\mu,\delta>0$
such that $\overline B_{R_0}(x_0)\subset \Omega$
and $B_{R_0}^\delta(x_0):=\{x\in B_{R_0}(x_0);\text{dist}(x,\partial B_{R_0}(x_0))<\delta\}$.
Here, $(R_0^2-|x-x_0|^2)_+=\max\{0,R_0^2-|x-x_0|^2\}$.
Then
$$\text{supp}\,u(x,t)=
\{(\theta,\rho(\theta,t));\theta\in S^{N-1}\},$$
where $(\theta,\rho)$ is the spherical coordinate centered at $x_0$,
$\rho(\theta,0)=R_0$ for all $\theta\in S^{N-1}$, and the propagating speed
$$
\pd{\rho(\theta,t)}{t}\Big|_{t=0}=R_0\Big(\frac{2m}{m-1}K_0^{m-1}-\chi\mu\Big),
\quad \forall \theta\in S^{N-1}.
$$
\end{theorem}

With the signal being consumed as time grows,
we show that the cells will eventually expand.

\begin{theorem}[Eventual expanding] \label{th-expand}
Let $(u,v)$ be a globally bounded weak solution of \eqref{eq-model}
with (i) $N=1$; or
(ii) $\sup_{t\in(0,\infty)}\|u(\cdot,t)\|_{C^{1/(2m)}(\overline\Omega)}\le C$ for some
constant $C>0$.
Further we assume the initial data $u_0\ge0$, $u_0\not\equiv0$
and $\Omega$ is convex.
Then there exist $\hat T>\hat t>0$ and $t_0\in(\hat t,\hat T)$, $\varepsilon_0>0$,
and a family of expanding open sets $\{A(t)\}_{t\in(\hat t,\hat T)}$, such that
$$A(t)\subset\text{supp}\,u(x,t),\quad t\in(\hat t,\hat T),$$
and $A(t)=\Omega$, $u(x,t)\ge\varepsilon_0$
for all $x\in\Omega$ and $t\in[t_0,\hat T]$.
\end{theorem}

Theorem \ref{th-expand} implies that the cells will eventually expand
to the whole domain.
After that we can show the eventual smoothness and large time behavior.

\begin{theorem}[Eventual smoothness] \label{th-smooth}
Let the assumptions in Theorem \ref{th-expand} be valid.
Then $u(x,t)\ge\varepsilon_0$ for all $x\in\Omega$ and $t\ge t_0$
with $t_0>0$ and $\varepsilon_0>0$ in Theorem \ref{th-expand},
$u\in C^{2,1}(\overline\Omega\times[t_0,\infty))$
and there exist $C>0$ and $c>0$ such that
$$\|u(\cdot,t)-\overline u\|_{L^\infty(\Omega)}
+\|v(\cdot,t)\|_{W^{1,\infty}(\Omega)}
\le Ce^{-ct}, \quad t>0,$$
where $\overline u=\int_\Omega u_0dx/|\Omega|$.
\end{theorem}

The main difficulty lies in the balance between the degenerate diffusion (expanding)
and the possible aggregating effect (shrinking) caused by the chemotaxis.
According to the exact propagating speed Theorem \ref{th-exact},
it is clear that the profile near the boundary of its support
competes with the gradient of the signal concentration.
We first prove the comparison principle by the approximate Hohmgren's approach,
and then construct several kinds of lower and upper solutions.
The self similar weak lower and upper solutions with
shrinking or expanding support are comparable with
the Barenblatt solution to the porous medium equation
\begin{equation} \label{eq-Barenblatt}
B(x,t)=(1+t)^{-k}\Big[\Big(1-\frac{k(m-1)}{2mN}\frac{|x|^2}{(1+t)^{2k/N}}
\Big)_+\Big]^\frac{1}{m-1}
\end{equation}
with $k=1/(m-1+2/N)$ for $m>1$.
After showing the eventual expanding property, we formulate
the eventual smoothness and large time behavior.

\section{Preliminaries: comparison principle and H\"older continuity}

\subsection{Comparison principle of degenerate diffusion equations}

We present the following comparison principle of degenerate diffusion equation in general form
\begin{equation} \label{eq-dege}
\begin{cases}
\displaystyle
\pd{u}{t}=\Delta A(u)-\nabla\cdot(B(u)\Phi(x,t)), \quad x\in\Omega, ~t>0,\\
(\nabla A(u)-B(u)\Phi)\cdot n=0, \quad x\in\partial\Omega, ~t>0,\\
u(x,0)=u_0(x), \quad x\in\Omega,
\end{cases}
\end{equation}
where $A(s)$ is strictly increasing and locally Lipchitz continuous for $s\in\mathbb R$,
$B(s)$ is locally Lipchitz continuous for $s\in\mathbb R$, and
$\Phi:\mathbb{R}^N\times\mathbb{R}_+\to\mathbb{R}^N$ is bounded.
Here the degenerate set $\{s\in\mathbb R;A'(s)=0\}$ has no interior point
and the equation \eqref{eq-dege} is weakly degenerate.
The typical case is $A(u)=u^m$ with $m>1$, $B(u)=\chi u$ and the solution $u$ is non-negative
(otherwise, one may write $A(u)=|u|^{m-1}u$).

\begin{lemma}[Comparison principle] \label{le-comparisonprinciple}
Let $T>0$ and the function space
$E=\{u\in L^\infty(Q_T);\nabla A(u)\in L^2(Q_T)\},$
$u_1,u_2\in E$, $\Phi\in L^\infty(Q_T)$,
and $u_1$, $u_2$ satisfy the following differential inequality
\begin{align*}
\begin{cases}
\displaystyle
\pd{u_1}{t}-\Delta A(u_1)+\nabla\cdot(B(u_1)\Phi(x,t)), &\\[2mm]
\displaystyle
\qquad\ge\pd{u_2}{t}-\Delta A(u_2)+\nabla\cdot(B(u_2)\Phi(x,t)),\qquad &x\in\Omega, t\in(0,T), \\[2mm]
\displaystyle
(\nabla A(u_1)-B(u_1)\Phi)\cdot n\ge(\nabla A(u_2)-B(u_2)\Phi)\cdot n, \qquad &x\in\partial\Omega, t\in(0,T),\\
u_1(x,0)\ge u_2(x,0), \qquad &x\in\Omega,
\end{cases}
\end{align*}
in the sense that the following inequality
\begin{align*}
&\iint_{Q_T} u_1\varphi_tdxdt
+\int_\Omega u_{10}(x)\varphi(x,0) dx
-\iint_{Q_T}  \nabla A(u_1)\cdot\nabla\varphi dxdt\\
&+\iint_{Q_T}  B(u_1)\Phi(x,t)\cdot\nabla\varphi dxdt
+\iint_{\partial\Omega\times(0,T)}(\nabla A(u_1)-B(u_1)\Phi)\cdot n dxdt, \\
\le&\iint_{Q_T}  u_2\varphi_tdxdt
+\int_\Omega u_{20}(x)\varphi(x,0) dx
-\iint_{Q_T}  \nabla A(u_2)\cdot\nabla\varphi dxd\\
&+\iint_{Q_T}  B(u_2)\Phi(x,t)\cdot\nabla\varphi dxdt
+\iint_{\partial\Omega\times(0,T)}(\nabla A(u_2)-B(u_2)\Phi)\cdot n dxdt,
\end{align*}
hold for some fixed $u_{10},u_{20}\in L^2(\Omega)$
such that $u_{10}\ge u_{20}$ on $\Omega$ and all test functions
$0\le\varphi\in L^2((0,T);W^{1,2}(\Omega))\cap W^{1,2}((0,T);L^2(\Omega))$ with
$\varphi(x,T)=0$ on $\Omega$.
Then $u_1(x,t)\ge u_2(x,t)$ almost everywhere in $Q_T$.
\end{lemma}
{\it \bfseries Proof.}
The following inequality
\begin{align*}
\iint_{Q_T}(u_1-u_2)\varphi_tdxdt
\le&\iint_{Q_T} \nabla (A(u_1)-A(u_2))\cdot\nabla\varphi dxdt\\
&-\iint_{Q_T} (B(u_1)-B(u_2))\Phi(x,t)\cdot\nabla\varphi dxdt,
\end{align*}
holds for all
$0\le\varphi\in L^2((0,T);W^{1,2}(\Omega))\cap W^{1,2}((0,T);L^2(\Omega))$ with
$\varphi(x,T)=0$.
If we further assume that $\pd{\varphi}{n}=0$ for $x\in\partial\Omega$ and $t\in(0,T)$,
then we have
\begin{equation} \label{eq-zcomp}
\iint_{Q_T}(u_1-u_2)\varphi_tdxdt
\le\iint_{Q_T} \Big(-(A(u_1)-A(u_2))\Delta\varphi
-(B(u_1)-B(u_2))\Phi(x,t)\cdot\nabla\varphi \Big)dxdt.
\end{equation}
Let
$$a(x,t)=\int_0^1 A'(su_1+(1-s)u_2)ds=
\begin{cases}
\displaystyle
\frac{A(u_1)-A(u_2)}{u_1-u_2}, ~&u_1(x,t)\ne u_2(x,t),\\
A'(u_1), ~&u_1(x,t)=u_2(x,t),
\end{cases}$$
\begin{align*}
b(x,t)&=\int_0^1 B'(su_1+(1-s)u_2)ds\cdot\Phi(x,t)\\
&\qquad=\begin{cases}
\displaystyle
\frac{(B(u_1)-B(u_2))\Phi(x,t)}{u_1-u_2},
~&u_1(x,t)\ne u_2(x,t),\\
B'(u_1)\Phi(x,t),~&u_1(x,t)=u_2(x,t),
\end{cases}
\end{align*}
and
$$c_\delta^\eta(x,t)=
\begin{cases}
(\eta+a(x,t))^{-\frac{1}{2}}b(x,t),\quad &|u_1(x,t)-u_2(x,t)|\ge\delta, \\
0,\quad &|u_1(x,t)-u_2(x,t)|<\delta,
\end{cases}$$
for any $\eta>0$ and $\delta>0$.
Further, for any fixed $\gamma>0$, we denote
$$F_\gamma=\{(x,t)\in Q_T;|u_1(x,t)-u_2(x,t)|\ge\gamma\},$$
and
$$G_\gamma=\{(x,t)\in Q_T;|u_1(x,t)-u_2(x,t)|<\gamma\}.$$
Now, \eqref{eq-zcomp} reads
\begin{equation} \label{eq-zcomp2}
\iint_{Q_T}(u_1-u_2)\Big(-\varphi_t
-a(x,t)\Delta\varphi
-b(x,t)\cdot\nabla\varphi \Big)dxdt\ge0,
\end{equation}
for all
$0\le\varphi\in L^2((0,T);W^{1,2}(\Omega))\cap W^{1,2}((0,T);L^2(\Omega))$ with
$\varphi(x,T)=0$ for $x\in\Omega$ and $\pd{\varphi}{n}=0$ for $x\in\partial\Omega$ and $t\in(0,T)$.
Since $\Phi(x,t),u_1,u_2$ are bounded and $A(s)$, $B(s)$ are locally Lipchitz continuous,
there exists a constant $C>0$ such that $|a|$, $|b|$ and $|u_1|$, $|u_2|\le C$.
Henceforth, a generic positive constant (possibly changing from line to line)
is denoted by $C$.
According to the strictly increasing property of $A(s)$
and the boundedness of $u_1$, $u_2$, there exists a constant $L(\gamma)>0$ such that
$$a(x,t)\ge L(\gamma), \quad \text{~for~all~}(x,t)\in F_\gamma,$$
and therefore
$$|c_\delta^\eta|\le L(\delta)^{-\frac{1}{2}}|b|\le L(\delta)^{-\frac{1}{2}}C=:K(\delta).$$

We employ the standard duality proof method
or the approximate Hohmgren's approach to complete this proof
(see Theorem 6.5 in \cite{PME}, Chapter 1.3 and 3.2 in \cite{NDE},
see also the comparison principle Lemma 3.4 in \cite{XuJDE} on unbounded domain
and Lemma 4.1 in \cite{XuMBE}).
For any smooth function $0\le\psi(x,t)\in C_0^2(Q_T)$,
consider the following approximated dual problem
\begin{equation} \label{eq-zdual}
\begin{cases}
-\varphi_t-(\eta+a_\varepsilon(x,t))\Delta\varphi
-c_{\delta,\varepsilon}^\eta(x,t)(\eta+a_\varepsilon(x,t))^\frac{1}{2}\cdot\nabla\varphi=\psi,
\quad (x,t)\in Q_T, \\
\displaystyle
\pd{\varphi}{n}=0, \quad (x,t)\in\partial\Omega\times(0,T),\\
\varphi(x,T)=0, \quad x\in\Omega,
\end{cases}
\end{equation}
where $\eta>0$, $\delta>0$, $\varepsilon>0$, $a_\varepsilon$ is a
smooth approximation of $a$ in $L^4(Q_T)$,
$a_\varepsilon\ge a$, and $c_{\delta,\varepsilon}^\eta(x,t)$ is a smooth approximation of
$c_{\delta}^\eta(x,t)$ in $L^4(Q_T)$.
Here we note that \eqref{eq-zdual} is a standard parabolic problem
as the initial data is imposed at the end time $t=T$.
Therefore, it has a smooth solution $\varphi\ge0$.
Maximum principle shows the boundedness of $\varphi$ such that $0\le\varphi\le C(\psi)$.
Then we get from \eqref{eq-zcomp2} and \eqref{eq-zdual} the estimate
\begin{align} \nonumber
\iint_{Q_T}(u_1-u_2)\psi dxdt\ge&
-\iint_{Q_T}|u_1-u_2||a-a_\varepsilon||\Delta\varphi|dxdt
\\ \nonumber
&-\eta\iint_{Q_T}|u_1-u_2||\Delta\varphi|dxdt
\\ \nonumber
&-\iint_{Q_T}|u_1-u_2||c_{\delta,\varepsilon}^\eta(\eta+a_\varepsilon)^\frac{1}{2}-b||\nabla\varphi| dxdt
\\ \label{eq-zI123}
=:&-I_1-I_2-I_3.
\end{align}

Next, we need the a priori estimate on $(\eta+a_\varepsilon)|\Delta\varphi|^2$.
We multiply the equation \eqref{eq-zdual} by $-\Delta\varphi$.
Integrating over $Q_T$ yields
\begin{align*}
&\iint_{Q_T}\varphi_t\Delta\varphi dxdt
+\iint_{Q_T}(\eta+a_\varepsilon)(\Delta\varphi)^2dxdt\\
\le& \iint_{Q_T} |c_{\delta,\varepsilon}^\eta|(\eta+a_\varepsilon)^\frac{1}{2}|\nabla\varphi||\Delta\varphi|dxdt
+\iint_{Q_T}\psi\Delta\varphi dxdt\\
\le&\frac{1}{4}\iint_{Q_T}(\eta+a_\varepsilon)(\Delta\varphi)^2dxdt+
\iint_{Q_T}|c_{\delta,\varepsilon}^\eta|^2|\nabla\varphi|^2dxdt
+\iint_{Q_T}|\Delta\psi||\varphi|dxdt\\
\le&\frac{1}{4}\iint_{Q_T}(\eta+a_\varepsilon)(\Delta\varphi)^2dxdt
+(K(\delta))^2\iint_{Q_T}|\nabla\varphi|^2dxdt
+C(\psi).
\end{align*}
Using $\varphi(x,T)=0$ and $\pd{\varphi}{n}=0$ on $\partial\Omega$, we have
\begin{align*}
\iint_{Q_T}\varphi_t\Delta\varphi dxdt
&=-\iint_{Q_T}\nabla\varphi\cdot\nabla\varphi_t dxdt
=-\frac{1}{2}\iint_{Q_T}\pd{}{t}|\nabla\varphi|^2 dxdt\\
&=\frac{1}{2}\int_\Omega |\nabla\varphi(x,0)|^2 dx\ge0,
\end{align*}
and
\begin{align*}
\iint_{Q_T}|\nabla\varphi|^2dxdt&=\iint_{Q_T}\nabla\varphi\cdot\nabla\varphi dxdt
=-\iint_{Q_T}\varphi\Delta\varphi dxdt\\
&\le \frac{1}{4(K(\delta))^2}\iint_{Q_T}(\eta+a_\varepsilon)(\Delta\varphi)^2dxdt
+\eta^{-1}(K(\delta))^2C(\psi).
\end{align*}
Therefore,
\begin{align} \label{eq-zDeltaphi}
(K(\delta))^2\iint_{Q_T}|\nabla\varphi|^2 dxdt
+\iint_{Q_T}(\eta+a_\varepsilon)(\Delta\varphi)^2dxdt
\le C(\psi)(K(\delta))^4\eta^{-1},
\end{align}
and
$$\|\Delta\varphi\|_{L^2(Q_T)}\le (C(\psi)(K(\delta))^4\eta^{-2})^\frac{1}{2}
\le C(\psi)(K(\delta))^2\eta^{-1}.$$
It follows that
\begin{align*}
I_1=&\iint_{Q_T}|u_1-u_2||a-a_\varepsilon||\Delta\varphi|dxdt\\
\le &C\|\Delta\varphi\|_{L^2(Q_T)}\|a-a_\varepsilon\|_{L^2(Q_T)}
\le C(\psi)(K(\delta))^2\eta^{-1}\|a-a_\varepsilon\|_{L^2(Q_T)},
\end{align*}
which converges to zero if we let $\varepsilon\to0$.
We can estimate $I_2$ as follows
\begin{align*}
I_2&=\eta\iint_{Q_T}|u_1-u_2||\Delta\varphi|dxdt
\\
&\le\eta\iint_{G_\gamma}|u_1-u_2||\Delta\varphi|dxdt
+\eta\iint_{F_\gamma}|u_1-u_2||\Delta\varphi|dxdt
\\
&\le\gamma\iint_{G_\gamma}\eta|\Delta\varphi|dxdt
+\frac{C\eta}{L(\gamma)^\frac{1}{2}}%
\iint_{F_\gamma}a^\frac{1}{2}|\Delta\varphi|dxdt
\\
&\le\gamma\iint_{G_\gamma}\eta|\Delta\varphi|dxdt
+\frac{C\eta}{L(\gamma)^\frac{1}{2}}%
\iint_{F_\gamma}a_\varepsilon^\frac{1}{2}|\Delta\varphi|dxdt
\\
&\le C\gamma\eta^{\frac{1}{2}}\Big(\iint_{Q_T}\eta|\Delta\varphi|^2dxdt\Big)^\frac{1}{2}
+\frac{C\eta}{L(\gamma)^\frac{1}{2}}%
\Big(\iint_{Q_T}a_\varepsilon|\Delta\varphi|^2dxdt\Big)^\frac{1}{2}
\\
&\le C\gamma\eta^{\frac{1}{2}}C(\psi)(K(\delta))^2\eta^{-\frac{1}{2}}
+\frac{C\eta}{L(\gamma)^\frac{1}{2}}C(\psi)(K(\delta))^2\eta^{-\frac{1}{2}}
\\
&=\gamma C(\psi)(K(\delta))^2+\eta^{\frac{1}{2}}C(\psi)(K(\delta))^2/L(\gamma)^\frac{1}{2}.
\end{align*}
We also have
\begin{align*}
I_3=&\iint_{Q_T}|u_1-u_2||c_{\delta,\varepsilon}^\eta(\eta+a_\varepsilon)^\frac{1}{2}-b||\nabla\varphi| dxdt
\\
\le& \iint_{G_\delta}|u_1-u_2||c_{\delta,\varepsilon}^\eta(\eta+a_\varepsilon)^\frac{1}{2}||\nabla\varphi| dxdt
+\iint_{G_\delta}|u_1-u_2||b||\nabla\varphi| dxdt
\\
&+\iint_{F_\delta}|u_1-u_2||c_{\delta,\varepsilon}^\eta(\eta+a_\varepsilon)^\frac{1}{2}-b||\nabla\varphi| dxdt
\\
\le& \delta\|c_{\delta,\varepsilon}^\eta(\eta+a_\varepsilon)^\frac{1}{2}\|_{L^2(G_\delta)}\|\nabla\varphi\|_{L^2(Q_T)}
+C\delta\iint_{G_\delta}|\nabla\varphi| dxdt
\\
&+C\|c_{\delta,\varepsilon}^\eta(\eta+a_\varepsilon)^\frac{1}{2}-b\|_{L^2(F_\delta)}\|\nabla\varphi\|_{L^2(Q_T)}.
\end{align*}
We note that
$$c_{\delta,\varepsilon}^\eta(\eta+a_\varepsilon)^\frac{1}{2}
\to c_{\delta}^\eta(\eta+a)^\frac{1}{2}
=\begin{cases}
0, \quad &(x,t)\in G_\delta, \\
b(x,t), \quad &(x,t)\in F_\delta,
\end{cases}$$
almost everywhere and also in $L^2(Q_T)$.
It follows that
$$
\limsup_{\varepsilon\to0}I_3\le C\delta\iint_{G_\delta}|\nabla\varphi| dxdt.
$$
We leave the uniform $L^1$ estimate of $\|\nabla\varphi\|_{L^1(Q_T)}\le C(\psi)$ to the next lemma
(Lemma \ref{le-L1}),
and we combine the above estimates to find
$$
\limsup_{\varepsilon\to0} (I_1+I_2+I_3)
\le \gamma C(\psi)(K(\delta))^2+\eta^{\frac{1}{2}}C(\psi)(K(\delta))^2/L(\gamma)^\frac{1}{2}
+C(\psi)\delta.
$$
Now we conclude according to \eqref{eq-zI123} that
\begin{align*}
\iint_{Q_T}(u_1-u_2)\psi dxdt\ge
-\Big\{\gamma C(\psi)(K(\delta))^2+\eta^{\frac{1}{2}}C(\psi)(K(\delta))^2/L(\gamma)^\frac{1}{2}
+C(\psi)\delta\Big\},
\end{align*}
for any given $\delta>0$, $\eta>0$, $\gamma>0$ and $\psi\ge0$,
which yields that
\begin{align*}
\iint_{Q_T}(u_1-u_2)\psi dxdt\ge0,
\end{align*}
by taking $\eta\to0$, then $\gamma\to0$, and at last $\delta\to0$.
Since $0\le\psi\in C_0^2(Q_T)$ is arbitrary selected, we see that $u_1\ge u_2$ almost everywhere on $Q_T$.
$\hfill\Box$

\begin{lemma} \label{le-L1}
Let $\varphi$ be the solution of the approximated dual problem \eqref{eq-zdual}
in the proof of Lemma \ref{le-comparisonprinciple}.
Then there holds
$$\sup_{t\in(0,T)}\int_{\Omega}|\nabla\varphi(x,t)| dx\le \iint_{Q_T}|\nabla \psi|dxdt.$$
\end{lemma}
{\it\bfseries Proof.}
Since $\varphi$ is smooth enough, $\varphi(x,T)=0$ on $\Omega$ and $\pd{\varphi}{n}=0$ on $\partial\Omega$,
we take the gradient of \eqref{eq-zdual} and then
multiply it by $|\nabla\varphi|^{\beta-1}\nabla \varphi$ with $\beta\in(0,1)$,
integrate over $Q_{t,T}=\Omega\times(t,T)$, to find
\begin{align}\nonumber
&\frac{1}{\beta+1}\int_\Omega|\nabla\varphi(x,t)|^{\beta+1}dx
+\beta\iint_{Q_{t,T}}(\eta+a_\varepsilon)|\Delta \varphi|^2|\nabla\varphi|^{\beta-1}dxdt
\\ \nonumber
=&-\beta\iint_{Q_{t,T}}c_{\delta,\varepsilon}^\eta(\eta+a_\varepsilon)^\frac{1}{2}\cdot\nabla\varphi
|\nabla\varphi|^{\beta-1}\Delta\varphi dxdt
+\iint_{Q_{t,T}}\nabla \psi\cdot|\nabla\varphi|^{\beta-1}\nabla \varphi dxdt
\\ \nonumber
\le& \beta\iint_{Q_{t,T}}(\eta+a_\varepsilon)|\Delta \varphi|^2|\nabla\varphi|^{\beta-1}dxdt
+\beta\iint_{Q_{t,T}}|c_{\delta,\varepsilon}^\eta|^2|\nabla\varphi|^{\beta+1}dxdt
\\ \label{eq-zL1}
&\qquad+\iint_{Q_{t,T}}|\nabla \psi||\nabla\varphi|^{\beta}dxdt.
\end{align}
According to \eqref{eq-zDeltaphi}, we see that
$$
\iint_{Q_T}|c_{\delta,\varepsilon}^\eta|^2|\nabla\varphi|^{\beta+1}dxdt
\le \iint_{Q_T}(K(\delta))^2(1+|\nabla\varphi|^{2})dxdt
\le C(\psi)(K(\delta))^4\eta^{-1},
$$
and
$$
\limsup_{\beta\to0}\iint_{Q_T}|\nabla \psi||\nabla\varphi|^{\beta}dxdt
\le\iint_{Q_T}|\nabla \psi|dxdt,
$$
by the dominated convergence theorem.
Now we let $\beta$ tends to zero, and \eqref{eq-zL1} implies that
$$
\int_\Omega|\nabla\varphi(x,t)|dx
\le \iint_{Q_T}|\nabla \psi|dxdt,
$$
for all $t\in(0,T)$.
The proof is completed.
$\hfill\Box$

The comparison principle together with specially constructed weak lower and upper solutions
are used to show the propagating properties.
Hence we define the following weak lower and upper solutions of the first equation in \eqref{eq-model}.

\begin{definition}[Weak lower and upper solutions]
A function $g(x,t)$ is said to be
a weak lower (or upper) solution of the first equation in \eqref{eq-model} on $Q_T$
corresponding to the initial value $u_0$ and
a given function $v$ such that $\nabla v\in L^\infty(Q_T)$,
if $0\le g\in L^\infty(Q_T)$, $\nabla g^m\in L^2(Q_T)$, and
it satisfies the following differential inequality
\begin{align*}
\begin{cases}
\displaystyle
\pd{g}{t}\le (\ge) \Delta g^m-\nabla\cdot(g\nabla v),\qquad &x\in\Omega, t\in(0,T), \\[3mm]
\displaystyle
\pd{g^m}{n}-g\nabla v\cdot n\le (\ge) 0, \qquad &x\in\partial\Omega, t\in(0,T),\\[2mm]
g(x,0)\ge0, \quad g(x,0)\le (\ge)u_0(x), \qquad &x\in\Omega,
\end{cases}
\end{align*}
where the first two inequality is satisfied
in the following sense
\begin{align*}
&\iint_{Q_T} g\varphi_tdxdt
+\int_\Omega g(x,0)\varphi(x,0) dx
\\
\ge (\le) &\iint_{Q_T} \nabla g^m\cdot\nabla\varphi dxdt
-\iint_{Q_T} g\nabla v\cdot\nabla\varphi dxdt,
\end{align*}
holds for all test functions
$0\le\varphi\in L^2((0,T);W^{1,2}(\Omega))\cap W^{1,2}((0,T);L^2(\Omega))$ with
$\varphi(x,T)=0$ on $\Omega$.
\end{definition}

\begin{lemma}[Comparison principle] \label{le-cp}
Let $(u,v)$ be a globally bounded weak solution of \eqref{eq-model}.
If $g(x,t)$ is a weak lower (or upper) solution of the first equation in \eqref{eq-model}
on $Q_T$, then
$$u(x,t)\ge (\le) g(x,t), \quad \forall (x,t)\in Q_T.$$
\end{lemma}
{\it\bfseries Proof.}
This is a simple corollary of comparison principle Lemma \ref{le-comparisonprinciple}.
$\hfill\Box$

\subsection{Regularity of H\"older continuity}

In order to show the propagation properties of the degenerate chemotaxis system \eqref{eq-model},
we need to know the existence, global boundedness, regularity and large time behavior
of its solutions.

We recall the existence and the global boundedness of solutions to the degenerate
chemotaxis model \eqref{eq-model}.

\begin{lemma}[\cite{JinarXiv}] \label{le-Jin}
Assume that $u_0\in L^\infty(\Omega)$, $\nabla u_0^m\in L^2(\Omega)$, $v_0\in W^{2,\infty}(\Omega)$,
$u_0,v_0\ge0$ and $m>1$, the spacial dimension $N=3$.
Then the problem \eqref{eq-model}
admits a nonnegative global bounded weak solution $(u,v)$ with
\begin{align*}
&\sup_{t\in(0,\infty)}(\|u(\cdot,t)\|_{L^\infty(\Omega)}+\|v\|_{W^{1,\infty}(\Omega)})\le C,
\\
&\sup_{t\in(0,\infty)}\int_\Omega|\nabla u^m|^2dx
+\sup_{t\in(0,\infty)}\|u^\frac{m+1}{2}\|_{W^{1,1}_2(\Omega\times(t,t+1))}\le C,
\\
&\sup_{t\in(0,\infty)}\|v\|_{W^{2,1}_p(\Omega\times(t,t+1))}\le C(p), \quad \forall p>1.
\end{align*}
Furthermore,
\begin{equation*}
\lim_{t\to\infty}\|v\|_{L^\infty(\Omega)}=0,
\quad \lim_{t\to\infty}\|u-\bar u\|_{L^p(\Omega)}=0, ~\forall p>1,
\end{equation*}
where $\bar u=\frac{1}{|\Omega|}\int_\Omega u_0 dx>0$.
\end{lemma}

\begin{remark}
The same global boundedness and asymptotic behavior results hold for the
lower spatial dimensional case $N=1,2$.
\end{remark}

\begin{remark}
We note that the boundedness of $\|u\|_{L^\infty(Q_T)}$ and $\|v\|_{L^\infty(Q_T)}$
is insufficient for the boundedness of $\|\Delta v\|_{L^\infty(Q_T)}$
according to the strong theory of the second equation in \eqref{eq-model}.
Hence the $W^{2,1}_p$ estimate for $p=\infty$
is not obtained in the above Lemma \ref{le-Jin}.
\end{remark}

\begin{remark}
One of the basic features for the degenerate diffusion equations,
such as the porous medium equation,
is the property of finite speed of propagation.
Therefore, the first component $u$ may not have positive minimum
for some time $t>0$.
For the large time behavior, it is proved in Lemma \ref{le-Jin} that
$u(x,t)$ converges to $\bar u$ in $L^p(\Omega)$ for $p<\infty$,
while the $L^\infty(\Omega)$ and some other more regular convergence are not deduced.
\end{remark}

In a special case that $v_0\equiv0$, we see that $v(x,t)\equiv0$ and
$u$ satisfies the porous medium equation.
The Barenblatt solution \eqref{eq-Barenblatt} of the porous medium equation
shows that the best regularity of the first equation in \eqref{eq-model} is no better than
H\"older continuous $C^\frac{1}{m-1}(\overline Q_T)$
(for $m>2$)
even for the one spatial dimensional case $N=1$.
In what follows, we will show the H\"older continuous of $u$ with respect to space,
and the boundedness of $\|\Delta v\|_{L^\infty(Q_T)}$.
Actually, we will prove that $\Delta v\in C^{\alpha,\alpha/2}(\overline Q_T)$
for some $\alpha\in(0,1)$.

\begin{lemma} \label{le-L4}
Let $N=1$ and $(u,v)$ be the globally bounded weak solution of \eqref{eq-model}.
Then there exists a constant $C>0$ such that
$$
\sup_{t\in(0,\infty)}\Big\{\|u^m(\cdot,t)\|_{C^{1/2}(\overline\Omega)}
+\|u(\cdot,t)\|_{C^{1/(2m)}(\overline\Omega)}\Big\}\le C.
$$
\end{lemma}
{\it\bfseries Proof.}
According to Lemma \ref{le-Jin}, $\|\nabla u^m(\cdot,t)\|_{L^2(\Omega)}$ is uniformly bounded.
The Sobolev embedding theorem for one dimensional case implies the uniform boundedness
of $\|u^m(\cdot,t)\|_{C^{1/2}(\overline\Omega)}$.

We assert that for $m>1$,
$$
|a-b|^m\le C(M)|a^m-b^m|, \quad \forall a,b\in [0,M].
$$
This is a simple result of calculus.
Actually, we can choose $C(M)=1$.
Therefore,
$$
\Big(\frac{|u(x_1,t)-u(x_2,t)|}{|x_1-x_2|^{1/(2m)}}\Big)^m
\le C(\sup_{t\in(0,\infty)}\|u\|_{L^\infty(\Omega)})
\frac{|u^m(x_1,t)-u^m(x_2,t)|}{|x_1-x_2|^{1/2}}, \quad x_1\ne x_2.
$$
That is, the uniform $C^{1/2}$ regularity of $u^m(\cdot,t)$ implies
the uniform $C^{1/(2m)}$ regularity of $u(\cdot,t)$.
$\hfill\Box$

The following continuity of $\|\nabla v(\cdot,t)\|_{L^\infty(\Omega)}$
and boundedness of $\|\Delta v(\cdot,t)\|_{L^\infty(\Omega)}$
will be used to formulate varies types of upper and lower solutions in the next section.

\begin{lemma} \label{le-Holder}
Let $(u,v)$ be the globally bounded weak solution of \eqref{eq-model}
such that
$$
\sup_{t\in(0,\infty)}\|u(\cdot,t)\|_{C^{1/(2m)}(\overline\Omega)}\le C,
$$
and $v_0\in C^{2,\alpha_0}(\overline\Omega)$ for some $\alpha_0\in(0,1)$,
$\pd{v_0}{n}=0$ on $\partial\Omega$.
Then $\nabla v(\cdot,t)$ is continuous
in the $\|\cdot\|_{L^\infty(\Omega)}$ norm with respect to time
and there exist $\alpha\in(0,1)$ and a constant $C(T,\delta)>0$ such that
$$
\|\Delta v(x,t)\|_{C^\alpha(\overline\Omega_\delta\times[0,T])}\le C(T,\delta),
$$
where $\Omega_\delta=\{x\in\Omega;\text{dist}(x,\partial\Omega)>\delta\}$.
\end{lemma}
{\it\bfseries Proof.}
Since $\|v\|_{W^{2,1}_p(\Omega\times(t,t+1))}$ is uniformly bounded for $p>1$
in Lemma \ref{le-Jin},
we see that $\sup_{t\in(0,\infty)}\|v(\cdot,t)\|_{C^\beta(\overline\Omega)}\le C$
for some $\beta\in(0,1)$.
Therefore,
$$\sup_{t\in(0,\infty)}\|(uv)(\cdot,t)\|_{C^\alpha(\overline\Omega)}\le C$$
for some $\alpha\in(0,1)$.
Indeed, we can choose $\alpha=\min\{1/(2m),\beta\}$.
The Schauder theory via Campanato space theory in \cite{Lieberman} implies
the interior H\"older continuity of $\Delta v$ with respect to space and time,
and the H\"older continuity of $v_t$ with respect to space
(the H\"older continuity of $v_t$ with respect to time is insufficient).
$\hfill\Box$

For large time behavior, we present the following regularity.

\begin{lemma} \label{le-large}
Let $(u,v)$ be the globally bounded weak solution of \eqref{eq-model}.
Then
$$\lim_{t\to\infty}\|\nabla v(\cdot,t)\|_{L^\infty(\Omega)}=0.$$
\end{lemma}
{\it\bfseries Proof.}
Let $(e^{t\Delta})_{t\ge0}$ be the Neumann heat semigroup in $\Omega$,
and let $\lambda_1>0$ denote the first nonzero eigenvalue of $-\Delta$
in $\Omega$ under Neumann boundary condition.
Then the solution $v$ can be expressed as follows
$$
v(x,t)=e^{t\Delta}v_0(x)-\int_0^t e^{(t-s)\Delta}(uv)(x,s)ds,
\quad t\ge t_0\ge0.
$$
According to the $L^p-L^q$ estimates for the Neumann heat semigroup
(see for example \cite{Winkler-Aggregation}),
\begin{align*}
\|\nabla v(x,t)&\|_{L^\infty(\Omega)}
\le \|\nabla e^{t\Delta}v_0(x)\|_{L^\infty(\Omega)}
+\int_0^t \|\nabla e^{(t-s)\Delta}(uv)(x,s)\|_{L^\infty(\Omega)}ds
\\
\le &C\big(1+t^{-\frac{1}{2}}\big)e^{-\lambda_1t}\|v_0\|_{L^\infty(\Omega)}
+\int_0^t C\big(1+(t-s)^{-\frac{1}{2}}\big)e^{-\lambda_1(t-s)}
\|(uv)(\cdot,s)\|_{L^\infty(\Omega)}ds
\\
\le &C\big(1+t^{-\frac{1}{2}}\big)e^{-\lambda_1t}\|v_0\|_{L^\infty(\Omega)}
+C\int_0^{t-1}e^{-\lambda_1(t-s)}ds
\\
&\qquad+C\int_{t-1}^t\big(1+(t-s)^{-\frac{1}{2}}\big)ds
\sup_{\tau\in(t-1,t)}\|v(\cdot,\tau)\|_{L^\infty(\Omega)},
\end{align*}
which tends to zero since $\|u(\cdot,t)\|_{L^\infty(\Omega)}$,
$\|v(\cdot,t)\|_{L^\infty(\Omega)}$ are uniformly bounded and
$\|v(\cdot,t)\|_{L^\infty(\Omega)}$ tends to zeros as $t\to\infty$
from Lemma \ref{le-Jin}.
$\hfill\Box$

\begin{lemma} \label{le-Deltav}
Let the conditions in Lemma \ref{le-Holder} be valid.
Then
$$\lim_{t\to\infty}\|\Delta v(\cdot,t)\|_{L^\infty(\Omega)}=0.$$
\end{lemma}
{\it\bfseries Proof.}
We rewrite $v=v_1+w$ such that
\begin{equation*}
\begin{cases}
v_{1t}=\Delta v_1, \quad x\in\Omega, ~t>0,\\
v_1(x,0)=v_0(x), \quad x\in\Omega, \\
\pd{v_1}{n}=0, \quad x\in\partial\Omega, ~t>0,
\end{cases}
\end{equation*}
and
\begin{equation*}
\begin{cases}
w_{t}=\Delta w-uv, \quad x\in\Omega, ~t>0,\\
w(x,0)=0, \quad x\in\Omega, \\
\pd{w}{n}=0, \quad x\in\partial\Omega, ~t>0.
\end{cases}
\end{equation*}
The Neumann heat semigroup theory shows
$\lim_{t\to\infty}\|\Delta v_1(\cdot,t)\|_{L^\infty(\Omega)}=0$.
We note that
$$\|(uv)(\cdot,t)\|_{C^\alpha(\overline\Omega)}
\le \|u(\cdot,t)\|_{L^\infty(\Omega)}\|v(\cdot,t)\|_{C^\alpha(\overline\Omega)}
+\|v(\cdot,t)\|_{L^\infty(\Omega)}\|u(\cdot,t)\|_{C^\alpha(\overline\Omega)}
\to0,$$
as $t$ tends to infinity since $\lim_{t\to\infty}\|v(\cdot,t)\|_{W^{1,\infty}(\Omega)}=0$
according to Lemma \ref{le-large}
and $\|u(\cdot,t)\|_{C^\alpha(\overline\Omega)}$ are uniformly bounded
in Lemma \ref{le-Holder} for some $\alpha\in(0,1)$.
The Schauder theory in \cite{Lieberman} shows the H\"older continuity
$$\|\Delta w(\cdot,t)\|_{L^\infty(\Omega)}
\le C_1\sup_{s\in[t/2,t]}\|(uv)(\cdot,s)\|_{C^\alpha(\overline\Omega)}
+C_2(t)\sup_{s\in(0,\infty)}\|(uv)(\cdot,s)\|_{C^\alpha(\overline\Omega)},$$
where $C_1>0$ is a constant and $C_2(t)$ decays to zeros as $t$ tends to infinity.
$\hfill\Box$

\section{Propagation properties: shrinking versus expanding}

This section is devoted to the study of the propagating properties of bacteria cells and
the large time behavior of the weak solution
$(u,v)$ to the problem \eqref{eq-model}.
In contrast with the heat equation,
it is known that the porous medium equation has the
property of finite speed of propagation.
Therefore, the first component $u$ may not have positive minimum
for some time $t>0$.
We use the comparison principle together with weak lower solutions.

Our interest lies in the propagating properties of the cell invasions.
Let us first focus on the waiting time and initial shrinking of the compact support.
Our approach is the combination of the comparison principle Lemma \ref{le-cp}
and weak lower and upper solutions with compact support.

\subsection{Initial shrinking caused by the chemotaxis}

The Barenblatt solution \eqref{eq-Barenblatt} of the classical porous medium equation
indicates the slow diffusion with finite speed of expanding support;
while the chemotaxis may cause backward diffusion, i.e. the aggregation,
which in competition with the slow diffusion results in a initial shrinking of the support
provided specified structures of the signal concentration.

We consider a typical situation in which the cells are concentrated in a compact support
and the signal concentration has the aggregation effect.
Specifically speaking, assume that
\begin{equation} \label{eq-shrinking}
\begin{cases}
\text{supp}\,u_0\subset\overline B_{R_0}(x_0)\subset\Omega,
\quad u_0\le K_0(R_0^2-|x-x_0|^2)^{d_0}, ~ x\in B_{R_0}(x_0),
\\
\nabla v_0\cdot (x-x_0) \le -\mu |x-x_0|^2, ~ x\in B_{R_0}(x_0),
\end{cases}
\end{equation}
for some $x_0\in\Omega$ and positive constants $d_0\ge 1/(m-1)$, and $R_0,K_0,\mu>0$.

We construct self similar upper and lower solution with compact support
to show the propagating property.
We note that for the degenerate porous medium type equation
and the self similar function of the form
$g=[(1-|x|^2)_+]^d$ with $md>1$,
we can check that $\nabla g^m$ is continuous and
$\Delta g^m\in L^q(\Omega)$ for some $q>1$.
This shows that the differential inequality for an upper (or lower) solution
only need to be valid almost everywhere,
without the possible Radon measures on the boundary of its support,
which is completely different from the uniform parabolic cases.

\begin{lemma} \label{le-upper}
Let the conditions in Lemma \ref{le-Holder} be valid
with the initial values satisfying \eqref{eq-shrinking}
and $\chi\mu>\frac{4m}{m-1}K_0^{m-1}\max\{1,R_0^{2((m-1)d_0-1)}\}$.
Define a function
$$g(x,t)=\varepsilon(\tau+t)^\sigma
\Big[\Big(\eta^2-\frac{|x-x_0|^2}{(\tau+t)^\beta}\Big)_+\Big]^d,
\quad x\in\Omega,~t\ge0,$$
where $d=1/(m-1)$, $\beta,\sigma\in\mathbb R$,
$\varepsilon>0$, $\eta>0$, $\tau>0$.
Then by appropriately selecting $\beta<0$, $\sigma>0,$
$\varepsilon$, $\eta$ and $\tau$,
the support of $g(x,t)$ is contained in $\Omega$
and shrinks for
$t\in(0,t_0)$ with some $t_0>0$ and
the function $g(x,t)$ is
an upper solution of the first equation in \eqref{eq-model}
on $\Omega\times(0,t_0)$
corresponding to $v(x,t)$ and the initial date $u_0$.
Therefore, $u(x,t)\le g(x,t)$ and there exist
a family of shrinking open sets
$\{A(t)\}_{t\in(0,t_0)}$ such that
$$\text{supp}\,u(\cdot,t)\subset\overline A(t)\subset\Omega, \quad t\in(0,t_0),$$
and $\partial A(t)$ has a finite derivative with respect to $t$.
\end{lemma}
{\it\bfseries Proof.}
For simplicity, we let
$$h(x,t)=\Big(\eta^2-\frac{|x-x_0|^2}{(\tau+t)^\beta}\Big)_+, \quad x\in\Omega, ~t\ge0,$$
and
$$A(t)=\Big\{x\in\Omega;\frac{|x-x_0|^2}{(\tau+t)^\beta}<\eta^2\Big\},
\quad t\ge0.$$
Without loss of generality, we may assume that $x_0=0$
and write $B_R=B_R(0)$.
Straightforward computation shows that
\begin{align*}
g_t=&\sigma\varepsilon(\tau+t)^{\sigma-1}h^d
+\varepsilon(\tau+t)^{\sigma}dh^{d-1}\frac{\beta|x|^2}{(\tau+t)^{\beta+1}},
\\
\nabla g=&-\varepsilon(\tau+t)^{\sigma}dh^{d-1}\frac{2x}{(\tau+t)^{\beta}},
\\
\nabla g^m=&-\varepsilon^m(\tau+t)^{m\sigma}mdh^{md-1}\frac{2x}{(\tau+t)^{\beta}},
\\
\Delta g^m=&\varepsilon^m(\tau+t)^{m\sigma}md(md-1)h^{md-2}%
\frac{4|x|^2}{(\tau+t)^{2\beta}}-\varepsilon^m(\tau+t)^{m\sigma}mdh^{md-1}%
\frac{2N}{(\tau+t)^{\beta}},
\end{align*}
for all $x\in A(t)$ and $t>0$.
According to the initial condition \eqref{eq-shrinking}
and the regularity result Lemma \ref{le-Holder},
we see that at the initial time
\begin{align*}
\nabla g(x,0)\cdot\nabla v(x,0)&=\nabla g(x,0)\cdot\nabla v_0(x)
\\
&=-\varepsilon\tau^{\sigma}dh^{d-1}\frac{2x}{\tau^{\beta}}\cdot\nabla v_0(x)
\ge \varepsilon\tau^{\sigma-\beta}dh^{d-1}{2\mu|x|^2},
\end{align*}
and there exists a $\hat t>0$ by the continuity such that
\begin{align*}
\nabla v(x,t)\cdot x&\le -\frac{\mu}{2}|x|^2,
\quad x\in B_{R_0}\backslash B_{R_0/2}, ~t\in[0,t_0],
\\
\nabla v(x,t)\cdot x&\le \frac{\mu}{2}R_0^2,
\quad x\in B_{R_0/2}, ~t\in[0,\hat t].
\end{align*}
Therefore,
\begin{align} \nonumber
\nabla g(x,t)\cdot\nabla v(x,t)
&=-\varepsilon(\tau+t)^{\sigma}dh^{d-1}\frac{2x}{(\tau+t)^{\beta}}\cdot\nabla v(x,t)
\\ \label{eq-znablav}
&\ge \varepsilon(\tau+t)^{\sigma-\beta}dh^{d-1}{\mu|x|^2},
\quad x\in B_{R_0}\backslash B_{R_0/2}, ~t\in[0,\hat t].
\end{align}

Let $\tau>0$ to be determined and
\begin{equation} \label{eq-zuppercondi}
\eta^2=\frac{R_0^2}{\tau^\beta},
\quad t_0=\min\{\tau,\hat t\}.
\end{equation}
According to the definition of $g(x,t)$, we see that
$A(0)=B_{R_0}(0)$, $\text{supp}\,u_0\subset\overline A(0)\subset\Omega$,
and $A(t)\subset B_{R_0}(0)\subset\Omega$ for $t\in[0,t_0]$.
Therefore, $\pd{g}{n}=0$ and $\pd{g^m}{n}=0$ on $\partial\Omega$
for all $t\in(0,t_0)$, and
\begin{align*}
g(x,0)&=\varepsilon\tau^\sigma
\Big[\Big(\eta^2-\frac{|x|^2}{\tau^\beta}\Big)_+\Big]^d
=\varepsilon\tau^\sigma\Big(\frac{R_0^2}{\tau^\beta}-\frac{|x|^2}{\tau^\beta}\Big)^d
\cdot1_{B_{R_0}(0)}
\\
&=\varepsilon\tau^{\sigma-d\beta}(R_0^2-|x|^2)^d
\cdot1_{B_{R_0}(0)}
\ge K_0(R_0^2-|x|^2)^{d_0}\cdot1_{B_{R_0}(0)}
\ge u_0(x),
\quad x\in\Omega,
\end{align*}
provided that
\begin{equation} \label{eq-zuppercondi1}
\varepsilon\tau^{\sigma-d\beta}\ge K_0\max\{1,R_0^{2(d_0-d)}\}.
\end{equation}
In order to find a weak upper solution $g$, we only need to check the following
differential inequality on $A(t)$
\begin{align} \label{eq-zupperweak}
\pd{g}{t}\ge\Delta g^m-\chi\nabla\cdot(g\nabla v)
=\Delta g^m-\chi\nabla g\cdot\nabla v-\chi g\Delta v, \quad x\in A(t), ~t\in(0,t_0).
\end{align}
We denote $C_1=\|\nabla v\|_{L^\infty(\Omega\times[0,1])}$
and $C_2=\|\Delta v\|_{L^\infty(\Omega\times[0,1])}$ for convenience,
since they are bounded according to Lemma \ref{le-Holder}.
A sufficient condition of inequality \eqref{eq-zupperweak} is
\begin{align} \nonumber
&\sigma\varepsilon(\tau+t)^{\sigma-1}h^d
+\varepsilon(\tau+t)^{\sigma}dh^{d-1}\frac{\beta|x|^2}{(\tau+t)^{\beta+1}}
\\ \nonumber
&\quad+\varepsilon^m(\tau+t)^{m\sigma}mdh^{md-1}\frac{2N}{(\tau+t)^{\beta}}
+\chi\nabla g\cdot\nabla v
\\ \label{eq-zuppersimilar}
\ge&\varepsilon^m(\tau+t)^{m\sigma}md(md-1)h^{md-2}%
\frac{4|x|^2}{(\tau+t)^{2\beta}}
+C_2\chi\varepsilon(\tau+t)^{\sigma}h^{d}.
\end{align}
for all $x\in A(t)$, $t\in(0,t_0)$.
As we have chosen $d=1/(m-1)$, we rewrite \eqref{eq-zuppersimilar} into
\begin{align} \nonumber
&\sigma\varepsilon(\tau+t)^{\sigma-1}h
+\frac{\varepsilon\beta}{m-1}(\tau+t)^{\sigma}\frac{|x|^2}{(\tau+t)^{\beta+1}}
\\ \nonumber
&\quad+2N\frac{m}{m-1}\varepsilon^m(\tau+t)^{m\sigma}\frac{h}{(\tau+t)^{\beta}}
+h^{1-d}\chi\nabla g\cdot\nabla v
\\ \label{eq-zupperA}
\ge&\frac{m}{(m-1)^2}\varepsilon^m(\tau+t)^{m\sigma}%
\frac{4|x|^2}{(\tau+t)^{2\beta}}
+C_2\chi\varepsilon(\tau+t)^{\sigma}h,
\end{align}
for all $x\in A(t)$, $t\in(0,t_0)$.
For simplicity, we denote \eqref{eq-zupperA} by $LHS\ge RHS$.

Now, we give sufficient conditions of \eqref{eq-zupperA} to be valid on
$B_{R_0/2}$ and $B_{R_0}\backslash B_{R_0/2}$ respectively
(Note that $A(t)\subset B_{R_0}(0)$ for $t\in(0,t_0)$ as $\beta<0$).
For $x\in (B_{R_0}\backslash B_{R_0/2})\cap A(t)$ and $t\in(t,t_0)$,
we have according to the estimate \eqref{eq-znablav} that
\begin{align} \nonumber
LHS-&RHS\ge
\sigma\varepsilon(\tau+t)^{\sigma-1}h
+\frac{\varepsilon\beta}{m-1}(\tau+t)^{\sigma}\frac{|x|^2}{(\tau+t)^{\beta+1}}
\\ \nonumber
&+2N\frac{m}{m-1}\varepsilon^m(\tau+t)^{m\sigma}\frac{h}{(\tau+t)^{\beta}}
+\chi \varepsilon(\tau+t)^{\sigma-\beta}d{\mu|x|^2}
\\ \nonumber
&-\frac{m}{(m-1)^2}\varepsilon^m(\tau+t)^{m\sigma}%
\frac{4|x|^2}{(\tau+t)^{2\beta}}
-C_2\chi\varepsilon(\tau+t)^{\sigma}h
\\ \nonumber
\ge& \Big(\sigma+2N\frac{m}{m-1}\varepsilon^{m-1}(\tau+t)^{(m-1)\sigma-\beta+1}
-C_2\chi(\tau+t)\Big)\varepsilon(\tau+t)^{\sigma-1}h
\\ \nonumber
&+\Big(d\chi\mu+\frac{\beta}{m-1}(\tau+t)^{-1}
-\frac{4m}{(m-1)^2}\varepsilon^{m-1}(\tau+t)^{(m-1)\sigma-\beta}\Big)%
\varepsilon(\tau+t)^{\sigma-\beta}{|x|^2}
\\ \nonumber
\ge& \Big(\sigma+2N\frac{m}{m-1}\varepsilon^{m-1}\tau^{(m-1)\sigma-\beta+1}
-2C_2\chi\tau\Big)\varepsilon(\tau+t)^{\sigma-1}h
\\ \label{eq-zLHSa}
&+\Big(d\chi\mu+\frac{\beta}{m-1}\tau^{-1}
-\frac{4m}{(m-1)^2}\varepsilon^{m-1}(2\tau)^{(m-1)\sigma-\beta}\Big)%
\varepsilon(\tau+t)^{\sigma-\beta}{|x|^2}.
\end{align}
For $x\in (B_{R_0/2})\cap A(t)$ and $t\in(t,t_0)$,
we also have
\begin{align} \nonumber
LHS-RHS\ge&
\Big(\sigma+2N\frac{m}{m-1}\varepsilon^{m-1}\tau^{(m-1)\sigma-\beta+1}
-2C_2\chi\tau\Big)\varepsilon(\tau+t)^{\sigma-1}h
\\ \nonumber
&+\Big(\frac{\beta}{m-1}\tau^{-1}
-\frac{4m}{(m-1)^2}\varepsilon^{m-1}(2\tau)^{(m-1)\sigma-\beta}\Big)%
\varepsilon(\tau+t)^{\sigma-\beta}{|x|^2}
\\ \label{eq-zLHSb}
&-d\chi\mu\varepsilon(\tau+t)^{\sigma-\beta}R_0^2.
\end{align}
Let $\beta\in[-2\ln(4/3)/\ln2,0)$, i.e. $2^{\beta/2}\in[3/4,1)$.
For $t\in(0,t_0)$, we see that
$$A(t)=B_{\eta(\tau+t)^{\beta/2}}(0)=B_{R_0(1+t/\tau)^{\beta/2}}(0)\supset
B_{2^{\beta/2}R_0}(0)\supset B_{3R_0/4}(0).$$
Further if $x\in (B_{R_0/2})\cap A(t)=B_{R_0/2}$ and $t\in(0,t_0)\subset(0,\tau)$,
$$h(x,t)=\Big(\eta^2-\frac{|x|^2}{(\tau+t)^\beta}\Big)_+
=\frac{R_0^2}{\tau^\beta}-\frac{|x|^2}{(\tau+t)^\beta}
\ge \frac{R_0^2}{\tau^\beta}-\frac{(R_0/2)^2}{(2\tau)^\beta}
\ge\frac{5}{9}\frac{R_0^2}{\tau^\beta}.$$
Then \eqref{eq-zLHSb} reads
\begin{align} \nonumber
&LHS-RHS\ge
\Big(\sigma+2N\frac{m}{m-1}\varepsilon^{m-1}\tau^{(m-1)\sigma-\beta+1}
-2C_2\chi\tau\Big)\varepsilon(\tau+t)^{\sigma-1}\frac{5}{9}\frac{R_0^2}{\tau^\beta}
\\ \nonumber
&\qquad+\Big(\frac{\beta}{m-1}\tau^{-1}
-\frac{4m}{(m-1)^2}\varepsilon^{m-1}(2\tau)^{(m-1)\sigma-\beta}-4d\chi\mu\Big)%
\varepsilon(\tau+t)^{\sigma-\beta}\frac{R_0^2}{4}
\\ \nonumber
\ge& \Big[
\Big(\sigma+2N\frac{m}{m-1}\varepsilon^{m-1}\tau^{(m-1)\sigma-\beta+1}
-2C_2\chi\tau\Big)\frac{5}{9\tau^\beta}
\\ \label{eq-zLHSc}
&+\Big(\frac{\beta}{m-1}\tau^{-1}
-\frac{4m}{(m-1)^2}\varepsilon^{m-1}(2\tau)^{(m-1)\sigma-\beta}-4d\chi\mu\Big)%
\frac{(\tau+t)^{1-\beta}}{4}
\Big]\varepsilon(\tau+t)^{\sigma-1}R_0^2.
\end{align}

Let $\varepsilon>0$, $\beta\in[-2\ln(4/3)/\ln2,0)$,
$\sigma>0$, $\tau>0$, $\eta>0$ and $t_0>0$ be chosen such that
\eqref{eq-zuppercondi}, \eqref{eq-zuppercondi1} are valid and
\begin{align} \label{eq-zuppercondi2}
\begin{cases}
\displaystyle
\sigma+2N\frac{m}{m-1}\varepsilon^{m-1}\tau^{(m-1)\sigma-\beta+1}
-2C_2\chi\tau\ge0,
\\
\displaystyle
d\chi\mu+\frac{\beta}{m-1}\tau^{-1}
-\frac{4m}{(m-1)^2}\varepsilon^{m-1}(2\tau)^{(m-1)\sigma-\beta}\ge0,
\\
\displaystyle
\Big(\sigma+2N\frac{m}{m-1}\varepsilon^{m-1}\tau^{(m-1)\sigma-\beta+1}
-2C_2\chi\tau\Big)\frac{5}{9\tau^\beta}
\\
\displaystyle
\qquad+\Big(\frac{\beta}{m-1}\tau^{-1}
-\frac{4m}{(m-1)^2}\varepsilon^{m-1}(2\tau)^{(m-1)\sigma-\beta}-4d\chi\mu\Big)%
\frac{(2\tau)^{1-\beta}}{4}
\ge0.
\end{cases}
\end{align}
We can fix $\tau=1$, $\eta$ and $t_0$ to be determined by \eqref{eq-zuppercondi},
$\varepsilon=K_0\max\{1,R_0^{2(d_0-d)}\}$ as \eqref{eq-zuppercondi1} is valid,
and $\beta<0$ with $|\beta|$ being sufficiently small such that
the second inequality in \eqref{eq-zuppercondi2} is true since
$\chi\mu>\frac{4m}{m-1}K_0^{m-1}\max\{1,R_0^{2((m-1)d_0-1)}\}$,
and at last we choose $\sigma>0$ to be sufficiently large such that
the first and the third inequalities are satisfied.
Now, \eqref{eq-zuppercondi2} is valid for those parameters.
Then according to the inequalities \eqref{eq-zLHSa}, \eqref{eq-zLHSb}, \eqref{eq-zLHSc},
we find that
$$LHS-RHS\ge0, \qquad x\in B_{R_0}\cap A(t)=A(t), ~t\in(0,t_0),$$
which yields \eqref{eq-zupperweak}, \eqref{eq-zupperA},
and then $g(x,t)$ is an upper solution.

The comparison principle Lemma \ref{le-cp} implies that
$u(x,t)\le g(x,t)$ for all $x\in\Omega$ and $t\in(0,t_0)$.
Thus,
$$\text{supp}\,u(\cdot,t)\subset\overline A(t)
=\{x\in\Omega;|x-x_0|^2<\eta^2(\tau+t)^\beta\},
\quad t\in(0,t_0),$$
and
$$\partial A(t)=\{x\in\Omega;|x-x_0|=\eta(\tau+t)^\frac{\beta}{2}\},
\quad t\in(0,t_0),$$
which has finite derivative with respect to $t$.
The family of sets $\{A(t)\}_{t\in(0,t_0)}$ is shrinking with respect to $t$
since $\beta<0$.
$\hfill\Box$

\begin{remark}
We compare the self similar weak upper solution $g(x,t)$
in the proof of Lemma \ref{le-upper} to the Barenblatt solution of porous medium equation
$$B(x,t)=(1+t)^{-k}\Big[\Big(1-\frac{k(m-1)}{2mN}\frac{|x|^2}{(1+t)^{2k/N}}
\Big)_+\Big]^\frac{1}{m-1},$$
with $k=1/(m-1+2/N)$.
The Barenblatt solution $B(x,t)$ is decaying at the rate $(1+t)^{-1/(m-1+2/N)}$
in $L^\infty(\mathbb R^N)$ and the support is expanding at the rate $(1+t)^{k/N}$.
Here, the upper solution
is increasing at the rate $(\tau+t)^{\sigma}$ and its support is
shrinking at the rate $(\tau+t)^{\beta/2}$.
The increasing of $g(x,t)$ makes it possible to be an upper solution,
which is crucial in the proof.
\end{remark}

\subsection{Finite speed propagating and the exact propagating speed}

We have proved that the compact support may shrink if the signal concentration
satisfies a special structure such as \eqref{eq-shrinking}.
Now, We will show the finite speed propagating property
without assuming the special structure on signal concentration.
Assume that
\begin{equation} \label{eq-finite}
\text{supp}\,u_0\subset\overline B_{R_0}(x_0)\subset\Omega,
\quad u_0\le K_0(R_0^2-|x-x_0|^2)^{d_0}, \quad x\in B_{R_0}(x_0),
\end{equation}
for some $x_0\in\Omega$ and positive constants $d_0\ge 1/(m-1)$ and $R_0,K_0>0$.

\begin{lemma} \label{le-upper-finite}
Let the conditions in Lemma \ref{le-Holder} be valid
with the initial values satisfying \eqref{eq-finite}.
Define a function
$$g(x,t)=\varepsilon(\tau+t)^\sigma
\Big[\Big(\eta^2-\frac{|x-x_0|^2}{(\tau+t)^\beta}\Big)_+\Big]^d,
\quad x\in\Omega,~t\ge0,$$
where $d=1/(m-1)$, $\beta,\sigma\in\mathbb R$,
$\varepsilon>0$, $\eta>0$, $\tau>0$.
Then by appropriately selecting $\beta>0$, $\sigma>0$
$\varepsilon$, $\eta$ and $\tau$,
the support of $g(x,t)$ is contained in $\Omega$ for
$t\in(0,t_0)$ with some $t_0>0$ and
the function $g(x,t)$ is
an upper solution of the first equation in \eqref{eq-model}
on $\Omega\times(0,t_0)$ corresponding to $v(x,t)$ and the initial data $u_0$.
Therefore, $u(x,t)\le g(x,t)$ and there exist
a family of open sets $\{A(t)\}_{t\in(0,t_0)}$ such that
$$\text{supp}\,u(\cdot,t)\subset\overline A(t)\subset\Omega, \quad t\in(0,t_0),$$
and $\partial A(t)$ has a finite derivative with respect to $t$.
\end{lemma}
{\it\bfseries Proof.}
This proof is similar to the proof of Lemma \ref{le-upper},
except there is no structure condition \eqref{eq-znablav}
and we need minor modifications.
We still define $h(x,t)$ and $A(t)$ as in the proof of Lemma \ref{le-upper}
and we assume $x_0=0$ for simplicity.
Let
\begin{equation} \label{eq-zuppercondiA}
\eta^2=\frac{R_0^2}{\tau^\beta},
\quad \varepsilon\tau^{\sigma-d\beta}\ge K_0\max\{1,R_0^{2(d_0-d)}\},
\end{equation}
and $C_1$, $C_2$ be defined as in the proof of Lemma \ref{le-upper}.
We need to check the differential inequality \eqref{eq-zupperweak}
(i.e. \eqref{eq-zuppersimilar}).
A sufficient condition of \eqref{eq-zuppersimilar} is
\begin{align} \nonumber
&\sigma\varepsilon(\tau+t)^{\sigma-1}h
+\frac{\varepsilon\beta}{m-1}(\tau+t)^{\sigma}\frac{|x|^2}{(\tau+t)^{\beta+1}}
\\ \nonumber
&\quad+2N\frac{m}{m-1}\varepsilon^m(\tau+t)^{m\sigma}\frac{h}{(\tau+t)^{\beta}}
-C_1\chi\varepsilon(\tau+t)^{\sigma-\beta}\frac{2|x|}{m-1}
\\ \label{eq-zupperB}
\ge&\frac{m}{(m-1)^2}\varepsilon^m(\tau+t)^{m\sigma}%
\frac{4|x|^2}{(\tau+t)^{2\beta}}
+C_2\chi\varepsilon(\tau+t)^{\sigma}h,
\end{align}
for all $x\in A(t)$, $t\in(0,t_0)$.
For simplicity, we denote \eqref{eq-zupperB} by $LHS\ge RHS$.

According to \eqref{eq-finite}, $\overline B_{R_0}(0)\subset\Omega$,
there exists a $R>R_0$ such that $\overline B_{R_0}(0)\subset B_R(0)\subset\subset\Omega$.
Let $\hat t>0$ depending on $\beta$ and $\tau$ such that
\begin{equation} \label{eq-zt0}
\Big(1+\frac{\hat t}{\tau}\Big)^\beta\le \frac{R^2}{R_0^2}.
\end{equation}
Let $t_0=\min\{\tau,\hat t\}$.
We see that for $t\in(0,t_0)$,
$$\text{supp}\,g(x,t)=\overline A(t)=\overline B_{\eta(\tau+t)^{\beta/2}}\cap\Omega
=\overline B_{R_0(1+t/\tau)^{\beta/2}}\cap\Omega
\subset\overline  B_R\cap\Omega=\overline B_R\subset \Omega.$$
Then $\pd{g}{n}=0$ and $\pd{g^m}{n}=0$ on $\partial\Omega$
for all $t\in(0,t_0)$.
For $x\in A(t)\backslash B_{R_0/2}$ and $t\in(0,t_0)$, we have
\begin{align} \nonumber
&LHS-RHS\ge
\sigma\varepsilon(\tau+t)^{\sigma-1}h
+\frac{\varepsilon\beta}{m-1}(\tau+t)^{\sigma}\frac{|x|^2}{(\tau+t)^{\beta+1}}
\\ \nonumber
&\quad+2N\frac{m}{m-1}\varepsilon^m(\tau+t)^{m\sigma}\frac{h}{(\tau+t)^{\beta}}
-C_1\chi\varepsilon(\tau+t)^{\sigma-\beta}\frac{4|x|^2}{(m-1)R_0}
\\ \nonumber
&\quad\quad-\frac{m}{(m-1)^2}\varepsilon^m(\tau+t)^{m\sigma}%
\frac{4|x|^2}{(\tau+t)^{2\beta}}
-C_2\chi\varepsilon(\tau+t)^{\sigma}h
\\ \nonumber
\ge& \Big(\sigma+2N\frac{m}{m-1}\varepsilon^{m-1}(\tau+t)^{(m-1)\sigma-\beta+1}
-C_2\chi(\tau+t)\Big)\varepsilon(\tau+t)^{\sigma-1}h
\\ \nonumber
&+\Big(\frac{\beta}{m-1}(\tau+t)^{-1}
-\frac{4m}{(m-1)^2}\varepsilon^{m-1}(\tau+t)^{(m-1)\sigma-\beta}
-\frac{4C_1\chi}{(m-1)R_0}\Big)%
\varepsilon(\tau+t)^{\sigma-\beta}{|x|^2}
\\ \nonumber
\ge& \Big(\sigma+2N\frac{m}{m-1}\varepsilon^{m-1}\tau^{(m-1)\sigma-\beta+1}
\min\{1,2^{(m-1)\sigma-\beta+1}\}
-2C_2\chi\tau\Big)\varepsilon(\tau+t)^{\sigma-1}h
\\ \nonumber
&\qquad+\Big(\frac{\beta}{m-1}(2\tau)^{-1}
-\frac{4m}{(m-1)^2}\varepsilon^{m-1}\tau^{(m-1)\sigma-\beta}
\max\{1,2^{(m-1)\sigma-\beta}\}
\\ \label{eq-zupperC}
&\qquad\qquad-\frac{4C_1\chi}{(m-1)R_0}\Big)%
\varepsilon(\tau+t)^{\sigma-\beta}{|x|^2}.
\end{align}
We note that $B_{R_0/2}\subset B_{R_0}\subset A(t)$ for $t\in(0,t_0)$ since $\beta>0$.
For $x\in (B_{R_0/2})\cap A(t)$ and $t\in(t,t_0)$,
we find that
$$
h(x,t)=\Big(\eta^2-\frac{|x|^2}{(\tau+t)^\beta}\Big)_+
=\frac{R_0^2}{\tau^\beta}-\frac{|x|^2}{(\tau+t)^\beta}
\ge \frac{R_0^2}{\tau^\beta}-\frac{(R_0/2)^2}{\tau^\beta}
\ge\frac{3}{4}\frac{R_0^2}{\tau^\beta},
$$
then we also have
\begin{align} \nonumber
&LHS-RHS\ge
\sigma\varepsilon(\tau+t)^{\sigma-1}h
+\frac{\varepsilon\beta}{m-1}(\tau+t)^{\sigma}\frac{|x|^2}{(\tau+t)^{\beta+1}}
\\ \nonumber
&\quad+2N\frac{m}{m-1}\varepsilon^m(\tau+t)^{m\sigma}\frac{h}{(\tau+t)^{\beta}}
-C_1\chi\varepsilon(\tau+t)^{\sigma-\beta}\frac{R_0}{m-1}
\\ \nonumber
&\quad\quad-\frac{m}{(m-1)^2}\varepsilon^m(\tau+t)^{m\sigma}%
\frac{4|x|^2}{(\tau+t)^{2\beta}}
-C_2\chi\varepsilon(\tau+t)^{\sigma}h
\\ \nonumber
\ge& \Big(\sigma+2N\frac{m}{m-1}\varepsilon^{m-1}(\tau+t)^{(m-1)\sigma-\beta+1}
-C_2\chi(\tau+t)\Big)\varepsilon(\tau+t)^{\sigma-1}h
\\ \nonumber
&+\Big(\frac{\beta}{m-1}(\tau+t)^{-1}
-\frac{4m}{(m-1)^2}\varepsilon^{m-1}(\tau+t)^{(m-1)\sigma-\beta}\Big)%
\varepsilon(\tau+t)^{\sigma-\beta}{|x|^2}
\\ \nonumber
&\qquad\qquad-C_1\chi\varepsilon(\tau+t)^{\sigma-\beta}\frac{R_0}{m-1}
\\ \nonumber
\ge& \Big(\sigma-2C_2\chi\tau\Big)\varepsilon(\tau+t)^{\sigma-1}%
\frac{3}{4}\frac{R_0^2}{\tau^\beta}
-C_1\chi\varepsilon(\tau+t)^{\sigma-\beta}\frac{R_0}{m-1}
\\ \label{eq-zupperD}
&+\Big(\frac{\beta}{m-1}(2\tau)^{-1}
-\frac{4m}{(m-1)^2}\varepsilon^{m-1}\tau^{(m-1)\sigma-\beta}
\max\{1,2^{(m-1)\sigma-\beta}\}\Big)%
\varepsilon(\tau+t)^{\sigma-\beta}{|x|^2},
\end{align}
provided that $\sigma\ge2C_2\chi\tau$.

Let $\tau=1$, $\eta=R_0$, $\varepsilon=K_0\max\{1,R_0^{2(d_0-d)}\}$,
and $\beta=(m-1)\sigma$ with $\sigma>0$ being sufficiently large such that
\begin{align*}
\begin{cases}
\displaystyle
\frac{\beta}{2(m-1)}
-\frac{4m}{(m-1)^2}\varepsilon^{m-1}
-\frac{4C_1\chi}{(m-1)R_0}\ge0,
\\[3mm]
\displaystyle
(\sigma-2C_2\chi)\frac{3R_0}{4\tau^\beta}\min\{1,2^{\beta-1}\}
-\frac{C_1\chi}{m-1}\ge0.
\end{cases}
\end{align*}
Then \eqref{eq-zupperD} tells us $LHS\ge RHS$ for all
$x\in A(t)$ and $t\in(0,t_0)$.
It follows that $g(x,t)$ is an upper solution.
The comparison principle Lemma \ref{le-cp} completes the proof.
$\hfill\Box$

Lemma \ref{eq-Barenblatt} implies the finite speed propagating property of
the degenerate diffusion equation.
We will present the exact propagating speed for a special structure initial data.

\begin{lemma}[Exact propagating speed]
Let the conditions in Lemma \ref{le-Holder} be valid
with the initial values satisfying
\begin{equation} \label{eq-speed}
\begin{cases}
u_0= K_0\big[(R_0^2-|x-x_0|^2)_+\big]^{d}, ~ x\in \Omega,
\\
\nabla v_0\cdot (x-x_0)= -\mu |x-x_0|^2, ~ x\in B_{R_0}^\delta(x_0),
\end{cases}
\end{equation}
for some $x_0\in\Omega$ and positive constants $d=1/(m-1)$, $R_0,K_0,\mu,\delta>0$
such that $\overline B_{R_0}(x_0)\subset \Omega$
and $B_{R_0}^\delta(x_0):=\{x\in B_{R_0}(x_0);\text{dist}(x,\partial B_{R_0}(x_0))<\delta\}$.
Then
$$\text{supp}\,u(x,t)=
\{(\theta,\rho(\theta,t));\theta\in S^{N-1}\},$$
where $(\theta,\rho)$ is the spherical coordinate centered at $x_0$,
$\rho(\theta,0)=R_0$ for all $\theta\in S^{N-1}$, and the propagating speed
$$
\pd{\rho(\theta,t)}{t}\Big|_{t=0}=R_0\Big(\frac{2m}{m-1}K_0^{m-1}-\chi\mu\Big),
\quad \forall\theta\in S^{N-1}.
$$
\end{lemma}
{\it\bfseries Proof.}
Define
$$g_\pm(x,t)=\varepsilon(\tau+t)^{\sigma_\pm}
\Big[\Big(\eta^2-\frac{|x-x_0|^2}{(\tau+t)^{\beta_\pm}}\Big)_+\Big]^d,
\quad x\in\Omega,~t\ge0,$$
with $\varepsilon=K_0$, $\tau=1$, $\eta=R_0$, $\sigma_\pm\in\mathbb R$,
$\beta_\pm\in\mathbb R$ are to be determined.
We have
$$g_\pm(x,0)=K_0\big[(R_0^2-|x-x_0|^2)_+\big]^{d}=u_0, \quad x\in\Omega,$$
and $\pd{g_\pm}{n}=0$, $\pd{g_\pm^m}{n}=0$
on $\partial\Omega$ at least for a small time interval
since $\overline B_{R_0}\subset\Omega$.
Here we only aim to find the exact propagating speed
and we only need to construct upper and lower solutions
on a small time interval.
We note that
\begin{align*}
\nabla g_\pm(x,0)\cdot \nabla v_0
&=-\varepsilon(\tau+t)^{\sigma_\pm-\beta}dh^{d-1}2(x-x_0)\cdot\nabla v_0\\
&=2\mu\varepsilon(\tau+t)^{\sigma_\pm-\beta}dh^{d-1}|x-x_0|^2,
\end{align*}
for $x\in B_{R_0}^\delta(x_0)$.
Let
$$\beta=\frac{4m}{m-1}K_0^{m-1}-2\chi\mu,$$
and $\beta_\pm$ approach $\beta$ from above and below.
Take $\sigma_+>0$ sufficiently large and $\sigma_-<0$ with $|\sigma_-|$
being sufficiently large, we can check
as in the proof of Lemma \ref{le-upper} and next Lemma \ref{le-lower} that
$g_\pm(x,t)$ are upper and lower solutions for a small time interval
$(0,T_\pm)$, where $T_\pm>0$ depend on $|\beta_\pm-\beta|$.
Here we omit the details.
Then the comparison principle Lemma \ref{le-cp} implies that
there exists $\{A_{\beta_\pm}(t)\}_{t\in(0,T_\pm)}$ such that
$$A_{\beta_\pm}(t)=B_{R_0(1+t)^{\beta_\pm/2}}(x_0), \quad t\in(0,T_\pm),$$
and
$$
A_{\beta_-}(t)\subset \text{supp}\,u(x,t) \subset\overline A_{\beta_+}(t),
\quad t\in(0,T_\pm).
$$
Therefore,
$$
\pd{\rho(\theta,t)}{t}\Big|_{t=0}
\in[R_0\beta_-/2,R_0\beta_+/2].
$$
Since $\beta_\pm$ approach $\beta$, we have
$\pd{\rho(\theta,0)}{t}=R_0\beta/2$.
$\hfill\Box$

\subsection{Eventual smoothness and expanding}

The large time behavior in Lemma \ref{le-Jin} and Lemma \ref{le-large}
shows that $\|v(\cdot,t)\|_{W^{1,\infty}(\Omega)}$ tends to zero as time grows.
This indicates that the chemotaxis effect decays and the support will expand to
the whole domain.
Now we construct a self similar weak lower solution with expanding support.

\begin{lemma} \label{le-lower}
Let the conditions in Lemma \ref{le-Holder} be valid
with the initial data $u_0\ge0$, $u_0\not\equiv0$
and $\Omega$ is convex.
Define a function
$$g(x,t)=\varepsilon(\tau+t)^\sigma
\Big[\Big(\eta^2-\frac{|x-x_0|^2}{(\tau+t)^\beta}\Big)_+\Big]^d,
\quad x\in\Omega,~t>-\tau,$$
where $d=1/(m-1)$, $\beta>0$, $\sigma<0$, $\varepsilon$, $\eta>0$, $\tau\in\mathbb R$ and
$x_0\in\Omega$.
Then by appropriately selecting $\beta$, $\varepsilon,\tau$, $\sigma$, $\eta$ and $x_0$,
the function $g(x,t)$ is
a weak lower solution of the first equation in \eqref{eq-model}
on $\Omega\times(\hat t,\hat T)$ corresponding to $v(x,t)$ and $u_0$
for some $\hat T>\hat t>0$.
Therefore, $u(x,t)\ge g(x,t)$ and there exist $t_0\in(\hat t,\hat T)$, $\varepsilon_0>0$,
and a family of expanding open sets $\{A(t)\}_{t\in(\hat t,\hat T)}$,
such that
$$A(t)\subset\text{supp}\,u(x,t),\quad t\in(\hat t,\hat T),$$
and $A(t)=\Omega$, $u(x,t)\ge\varepsilon_0$
for all $x\in\Omega$ and $t\in[t_0,\hat T]$.
\end{lemma}
{\it\bfseries Proof.}
Since $u_0\ge0$, $u_0\not\equiv0$ and $u_0\in C(\overline\Omega)$,
the first equation in \eqref{eq-model} shows that
$$
\int_\Omega u(x,t)dx=\int_\Omega u_0(x)>0, \quad t>0.
$$
For any $t>0$, there exists a $x_0(t)\in\Omega$ such that
$u(x_0(t),t)\ge \overline u:=\frac{1}{|\Omega|}\int_\Omega u_0(x)>0$.
According to the uniform H\"older continuity of $u(\cdot,t)$,
we find that there exists a $R_0>0$ independent of $t$ such that
\begin{equation} \label{eq-zloweru0}
u(x,t)\ge \frac{\overline u}{2}=:\varepsilon_1, \quad \forall x\in B_{R_0}(x_0(t)).
\end{equation}
We denote $C_1(t)=\|\nabla v(\cdot,t)\|_{L^\infty(\Omega)}$
and $C_2(t)=\|\Delta v(\cdot,t)\|_{L^\infty(\Omega)}$ for convenience.
According to Lemma \ref{le-large} and Lemma \ref{le-Deltav},
$C_1(t)$ and $C_2(t)$ tend to zero.
For fixed $\delta>0$ to be determined, let $\hat t>0$ depend on $\delta$
such that
\begin{equation} \label{eq-zC1}
C_1(t)\le \delta, \quad C_2(t)\le \delta, \qquad \forall t\ge \hat t.
\end{equation}
Note that $u(x,\hat t)\ge\varepsilon_1$ on $B_{R_0}(x_0(\hat t))$.
Without loss of generality, we may assume that $B_{R_0}=B_{R_0}(x_0(\hat t))\subset\Omega$
and $x_0=x_0(\hat t)=0$.

Similar to the proof of Lemma \ref{le-upper}, we let
$$h(x,t)=\Big(\eta^2-\frac{|x-x_0|^2}{(\tau+t)^\beta}\Big)_+, \quad x\in\Omega, ~t\ge0,$$
and
$$A(t)=\Big\{x\in\Omega;\frac{|x-x_0|^2}{(\tau+t)^\beta}<\eta^2\Big\},
\quad t\ge0.$$
According to the definition of $g$, we see that
$\pd{g}{n}\le0$ and $\pd{g^m}{n}\le0$ on $\partial\Omega$
since $\Omega$ is convex, and for $\tau=1-\hat t$ we have
$$g(x,\hat t)=\varepsilon[(\eta^2-|x|^2)_+]^d\le \varepsilon_11_{B_{R_0}(x_0)}\le u_0(x),
\quad x\in\Omega,$$
provided that
\begin{equation} \label{eq-zcondi1}
\eta\le R_0, \quad \varepsilon\eta^{2d}\le\varepsilon_1.
\end{equation}
In order to find a weak lower solution $g$, we only need to check the following
differential inequality on $A(t)$
\begin{align} \label{eq-zweak}
\pd{g}{t}\le\Delta g^m-\chi\nabla\cdot(g\nabla v)
=\Delta g^m-\chi\nabla g\cdot\nabla v-\chi g\Delta v, \quad x\in A(t), ~t\in(\hat t,\hat T),
\end{align}
for some $\hat T>\hat t$ to be determined.

A sufficient condition of inequality \eqref{eq-zweak} is
\begin{align} \nonumber
&\sigma\varepsilon(\tau+t)^{\sigma-1}h
+\frac{\varepsilon\beta}{m-1}(\tau+t)^{\sigma}\frac{|x|^2}{(\tau+t)^{\beta+1}}
\\ \nonumber
&\quad+2N\frac{m}{m-1}\varepsilon^m(\tau+t)^{m\sigma}\frac{h}{(\tau+t)^{\beta}}
+C_1(t)\chi\varepsilon(\tau+t)^{\sigma-\beta}\frac{2|x|}{m-1}
\\ \label{eq-zupperE}
\le&\frac{m}{(m-1)^2}\varepsilon^m(\tau+t)^{m\sigma}%
\frac{4|x|^2}{(\tau+t)^{2\beta}}
-C_2(t)\chi\varepsilon(\tau+t)^{\sigma}h,
\end{align}
for all $x\in A(t)$, $t\in(\hat t,\hat T)$.
For simplicity, we denote \eqref{eq-zupperE} by $LHS\le RHS$.
The estimates on the above inequality is quit similar to \eqref{eq-zupperB}
in the proof of Lemma \ref{le-upper-finite}
except some terms are with inverse signs.
Here, \eqref{eq-zupperC} and \eqref{eq-zupperD} are changed into
\begin{align} \nonumber
&LHS-RHS\le
\sigma\varepsilon(\tau+t)^{\sigma-1}h
+\frac{\varepsilon\beta}{m-1}(\tau+t)^{\sigma}\frac{|x|^2}{(\tau+t)^{\beta+1}}
\\ \nonumber
&\quad+2N\frac{m}{m-1}\varepsilon^m(\tau+t)^{m\sigma}\frac{h}{(\tau+t)^{\beta}}
+C_1(t)\chi\varepsilon(\tau+t)^{\sigma-\beta}\frac{4|x|^2}{(m-1)R_0}
\\ \nonumber
&\quad\quad-\frac{m}{(m-1)^2}\varepsilon^m(\tau+t)^{m\sigma}%
\frac{4|x|^2}{(\tau+t)^{2\beta}}
+C_2(t)\chi\varepsilon(\tau+t)^{\sigma}h
\\ \nonumber
\le& \Big(\sigma+2N\frac{m}{m-1}\varepsilon^{m-1}(\tau+t)^{(m-1)\sigma-\beta+1}
+C_2(t)\chi(\tau+t)\Big)\varepsilon(\tau+t)^{\sigma-1}h
\\ \nonumber
&+\Big(\frac{\beta}{m-1}
-\frac{4m}{(m-1)^2}\varepsilon^{m-1}(\tau+t)^{(m-1)\sigma-\beta+1}
+\frac{4C_1(t)\chi(\tau+t)}{(m-1)R_0}\Big)%
\varepsilon(\tau+t)^{\sigma-\beta-1}{|x|^2}
\\ \nonumber
\le& \Big(\sigma+2N\frac{m}{m-1}\varepsilon^{m-1}
\max\{1,(\tau+\hat T)^{(m-1)\sigma-\beta+1}\}
+C_2(t)\chi(\tau+\hat T)\Big)\varepsilon(\tau+t)^{\sigma-1}h
\\ \nonumber
&\qquad+\Big(\frac{\beta}{m-1}
-\frac{4m}{(m-1)^2}\varepsilon^{m-1}
\min\{1,(\tau+\hat T)^{(m-1)\sigma-\beta+1}\}
\\ \label{eq-zupperC2}
&\qquad\qquad+\frac{4C_1(t)\chi(\tau+\hat T)}{(m-1)R_0}\Big)%
\varepsilon(\tau+t)^{\sigma-\beta-1}{|x|^2},
\quad x\in A(t)\backslash B_{R_0/2},~t\in(\hat t,\hat T),
\end{align}
and (note that $\sigma<0$)
\begin{align} \nonumber
&LHS-RHS\le
\sigma\varepsilon(\tau+t)^{\sigma-1}h
+\frac{\varepsilon\beta}{m-1}(\tau+t)^{\sigma}\frac{|x|^2}{(\tau+t)^{\beta+1}}
\\ \nonumber
&\quad+2N\frac{m}{m-1}\varepsilon^m(\tau+t)^{m\sigma}\frac{h}{(\tau+t)^{\beta}}
+C_1(t)\chi\varepsilon(\tau+t)^{\sigma-\beta}\frac{R_0}{m-1}
\\ \nonumber
&\quad\quad-\frac{m}{(m-1)^2}\varepsilon^m(\tau+t)^{m\sigma}%
\frac{4|x|^2}{(\tau+t)^{2\beta}}
+C_2(t)\chi\varepsilon(\tau+t)^{\sigma}h
\\ \nonumber
\le& \Big(\sigma+2N\frac{m}{m-1}\varepsilon^{m-1}(\tau+t)^{(m-1)\sigma-\beta+1}
+C_2(t)\chi(\tau+t)\Big)\varepsilon(\tau+t)^{\sigma-1}h
\\ \nonumber
&+\Big(\frac{\beta}{m-1}
-\frac{4m}{(m-1)^2}\varepsilon^{m-1}(\tau+t)^{(m-1)\sigma-\beta+1}\Big)%
\varepsilon(\tau+t)^{\sigma-\beta-1}{|x|^2}
\\ \nonumber
&\qquad\qquad+C_1(t)\chi\varepsilon(\tau+t)^{\sigma-\beta}\frac{R_0}{m-1}
\\ \nonumber
\le& \Big(2N\frac{m}{m-1}\varepsilon^{m-1}
\max\{1,(\tau+\hat T)^{(m-1)\sigma-\beta+1}\}
+C_2(t)\chi(\tau+\hat T)
\Big)\varepsilon(\tau+t)^{\sigma-1}\eta^2
\\ \nonumber
&+\Big(\frac{\beta}{m-1}
-\frac{4m}{(m-1)^2}\varepsilon^{m-1}
\min\{1,(\tau+\hat T)^{(m-1)\sigma-\beta+1}\}\Big)%
\varepsilon(\tau+t)^{\sigma-\beta-1}{|x|^2}
\\ \label{eq-zupperD2}
& +\sigma\varepsilon(\tau+t)^{\sigma-1}\frac{3}{4}\eta^2
+C_1(t)\chi\varepsilon(\tau+t)^{\sigma-\beta}\frac{R_0}{m-1},
\quad x\in (B_{R_0/2})\cap A(t), ~t\in(\hat t,\hat T).
\end{align}

Since $\Omega$ is bounded, there exists $R>R_0$ such that $\Omega\subset B_R(x_0)$.
Let $\eta=R_0$, $\varepsilon>0$, $\beta\in(0,1)$, $\tau=1-\hat t$, $\hat T>\hat t$
and $\sigma=-\frac{1-\beta}{m-1}<0$ be chosen such that
\begin{align} \label{eq-zcondi2}
\begin{cases}
\varepsilon\eta^{2d}\le\varepsilon_1, \quad
2N\frac{m}{m-1}\varepsilon^{m-1}\le -\sigma/4,
\quad \beta\le \frac{2m}{m-1}\varepsilon^{m-1},
\\
\delta\chi(\hat T-\hat t+1)\le -\sigma/4, \quad
4\delta\chi(\hat T-\hat t+1)\le \frac{2m}{m-1}\varepsilon^{m-1}R_0, \\
\delta\chi(\hat T-\hat t+1)^{1-\beta}\frac{R_0}{m-1}\le -\sigma\eta^2/4,
\quad (\hat T-\hat t+1)^{\beta/2}\ge 2R/R_0.
\end{cases}
\end{align}
The above seven inequalities can be satisfied simultaneously in the following way.
We first fix $\beta\in(0,1)$ sufficiently small such that
$N\beta\le (1-\beta)/(4(m-1))$.
Then we set $\varepsilon=\varepsilon(\beta)>0$
such that $\frac{2m}{m-1}\varepsilon^{m-1}=\beta$.
Now we can modify $\beta$ to be smaller such that
$\varepsilon\eta^{2d}\le\varepsilon_1$.
The first three inequalities are valid.
Let $L=e^{\frac{2}{\beta}\ln\frac{2R}{R_0}}-1$ and
$$
\delta=\min\{-\sigma/(4\chi(L+1)),\frac{2m}{m-1}\varepsilon^{m-1}R_0/(4\chi(L+1)),
-\sigma\eta^2(m-1)/(4\chi(L+1)^{1-\beta}R_0)\}.
$$
For this $\delta>0$, let $\hat t$ be chosen such that \eqref{eq-zC1} is fulfilled
and $\hat T=\hat t+L$.

For those parameters, we see that \eqref{eq-zcondi2} is valid
and \eqref{eq-zupperC2}, \eqref{eq-zupperD2} tells us
$LHS\le RHS$ for all $x\in A(t)$ and $t\in(\hat t,\hat T)$,
i.e. \eqref{eq-zupperE}.
It follows that $g(x,t)$ is a lower solution.
The comparison principle Lemma \ref{le-cp} shows that
$$u(x,t)\ge g(x,t)=\varepsilon(\tau+t)^\sigma
\Big[\Big(\eta^2-\frac{|x-x_0|^2}{(\tau+t)^\beta}\Big)_+\Big]^d,
$$
for all $x\in\Omega$ and $t\in(\hat t,\hat T)$.
We note that for this lower solution, its support satisfies
$$A(\hat t)=B_{\eta(\tau+\hat t)^{\beta/2}}(x_0)\cap\Omega=B_{R_0}(x_0),$$
and
$$A(\hat T)=B_{\eta(\tau+\hat T)^{\beta/2}}(x_0)\cap\Omega
=B_{R_0(\hat T-\hat t+1)^{\beta/2}}(x_0)\cap\Omega
\supset B_{2R}(x_0)\cap\Omega=\Omega,$$
since $(\hat T-\hat t+1)^{\beta/2}\ge 2R/R_0$ in \eqref{eq-zcondi2}
and $\Omega\subset B_R(x_0)$.
There exists a $t_1\in(\hat t,\hat T)$ such that
$$\eta^2-\frac{|x-x_0|^2}{(\tau+t)^\beta}\ge0,
\quad \forall x\in\Omega,~t\in(t_1,\hat T),$$
which means $A(t)=\Omega$ for $t\in(t_1,\hat T)$.
And there exists a $t_0\in(t_1,\hat T)$ such that
$$\eta^2-\frac{|x-x_0|^2}{(\tau+t)^\beta}\ge \frac{\eta^2}{2},
\quad \forall x\in\Omega,~t\in(t_0,\hat T),$$
and thus
$$u(x,t)\ge g(x,t)\ge\varepsilon(\hat T-\hat t+1)^{\sigma}
\Big(\frac{\eta^2}{2}\Big)^d=:\varepsilon_0,
\quad \forall x\in\Omega,~t\in(t_0,\hat T).$$
The proof is completed.
$\hfill\Box$

\begin{remark}
It is interesting to compare the self similar weak lower solution $g(x,t)$
in the proof of Lemma \ref{le-lower} to the Barenblatt solution of porous medium equation
$$B(x,t)=(1+t)^{-k}\Big[\Big(1-\frac{k(m-1)}{2mN}\frac{|x|^2}{(1+t)^{2k/N}}
\Big)_+\Big]^\frac{1}{m-1},$$
with $k=1/(m-1+2/N)$.
The Barenblatt solution $B(x,t)$ is decaying at the rate $(1+t)^{-1/(m-1+2/N)}$
in $L^\infty(\mathbb R^N)$ and the support is expanding at the rate $(1+t)^{k/N}$.
While the self similar weak lower solution $g(x,t)$
is decaying at the rate $(1+t)^{-(1-\beta)/(m-1)}$ and its support is
expanding at the rate $(1+t)^{\beta/2}$.
Here in the proof we have selected $\beta>0$ sufficiently small,
which means the support of $g$ is expanding with a much slower rate
and the maximum of $g$ is decaying at a slightly faster rate.
\end{remark}

Now that we have proved the lower bound of $u(x,t)$ on $\Omega\times(t_0,\hat T)$,
we will show the globally lower bound at large time,
as well as the non-degeneracy, regularity for large time behavior.

\begin{lemma}[Eventual smoothness]
Let the conditions in Lemma \ref{le-lower} be valid.
Then $u(x,t)\ge\varepsilon_0$ for all $x\in\Omega$ and $t\ge t_0$
with $t_0>0$ and $\varepsilon_0>0$ being defined as in the proof of Lemma \ref{le-lower},
$u\in C^{2,1}(\overline\Omega\times[t_0,\infty))$
and there exist $C>0$ and $c>0$ such that
$$\|u(\cdot,t)-\overline u\|_{L^\infty(\Omega)}
+\|v(\cdot,t)\|_{W^{1,\infty}(\Omega)}
\le Ce^{-ct}, \quad t>0,$$
where $\overline u=\int_\Omega u_0dx/|\Omega|$.
\end{lemma}
{\it\bfseries Proof.}
We point out that
$$\varepsilon_0=\varepsilon(\hat T-\hat t+1)^{\sigma}
\Big(\frac{\eta^2}{2}\Big)^d
=\varepsilon(L+1)^{\sigma}
\Big(\frac{\eta^2}{2}\Big)^d$$
is independent of $\delta$ and $\hat t$ therein,
since $L$ only depends on $\beta$, $R_0$ and $R$
(note that $\beta$, $\sigma$, $\varepsilon$ depend only on $\varepsilon_1$
and $\varepsilon_1=\overline u/2$ is fixed).
Therefore, we can take $\hat t$ larger to be $\hat t+\theta$ with any $\theta>0$
such that \eqref{eq-zC1} is also valid.
Lemma \ref{le-lower} shows that
$u(x,t)\ge\varepsilon_0$ for all $x\in\Omega$ and $t\in[t_0+\theta,\hat T+\theta]$.
Since $\varepsilon_0>0$ is fixed and $\theta>0$ is arbitrary,
we have $u(x,t)\ge\varepsilon_0$ for all $x\in\Omega$ and $t\ge t_0$.
It follows that the first equation in \eqref{eq-model} is non-degenerate and
uniform parabolic.
The H\"older regularity and exponential decay
can be verified similar to the proof of Theorem 1.3 in \cite{XuMBE}.
$\hfill\Box$

\section*{Acknowledgement}
The research of S. Ji is supported by NSFC Grant No. 11701184.
The research of M. Mei was supported in part
by NSERC Grant RGPIN 354724-2016, and FRQNT Grant No. 2019-CO-256440.
The research of J. Yin was supported in
part by NSFC Grant No. 11771156.

\end{document}